\documentclass{gtart_h}

\def\ifplaintex{\expandafter\ifx\csname documentclass\endcsname\relax}

\ifplaintex 
\else
\fi

\def\gt{{\mathsurround=0pt\it $\cal G\mskip-2mu$eometry \&\ 
$\cal T\!\!$opology}}        

\def\gtm{{\mathsurround=0pt\it $\cal G\mskip-2mu$eometry \&\ 
$\cal T\!\!$opology $\cal M\mskip-1mu$onographs}}    

\def\gtp{{\mathsurround=0pt\it $\cal G\mskip-2mu$eometry \&\ 
$\cal T\!\!$opology $\cal P\!$ublications}}  

\def\recd{{\small Received:\qua\receiveddate\ifx\reviseddate\relax
\else\qquad Revised:\qua\reviseddate\fi\par}} 

\def\volumenumber#1{\def\thevolumenumber{#1}}
\def\volumeyear#1{\def\thevolumeyear{#1}}
\def\volumename#1{\def\thevolumename{#1}}
\def\papernumber#1{\def\thepapernumber{#1}}
\def\pagenumbers#1#2{\def\startpage{#1}\def\finishpage{#2}}
\def\published#1{\def\publishdate{#1}}
\def\received#1{\def\receiveddate{#1}}
\def\revised#1{\def\reviseddate{#1}}
\def\accepted#1{\def\accepteddate{#1}}

\def\asciiaddress#1{\def\theasciiaddress{#1}}
\def\asciiemail#1{\def\theasciiemail{#1}}

\long\def\asciiabstract#1{\long\def\theasciiabstract{#1}}

\let\\\par
\let\thevolumenumber\relax\let\thepapernumber\relax
\let\thevolumeyear\relax\let\startpage\relax
\let\finishpage\relax\let\publishdate\relax\let\receiveddate\relax
\let\reviseddate\relax\let\accepteddate\relax\let\theasciititle\relax
\let\theasciiauthors\relax\let\theasciiaddress\relax
\let\theasciiabstract\relax

\let\theerratum\relax\let\theasciiemail\relax
\let\theshortauthors\relax\let\theshorttitle\relax

\def\startpage{1}\def\finishpage{15}\def\thepapernumber{77}

\volumenumber{2}
\volumename{Proceedings of the Kirbyfest}
\volumeyear{1999}

\long\def\maketitlep{   

\count0=\startpage

\gtm\nl        
{\small Volume \thevolumenumber: \thevolumename\nl 
\ifx\theerratum\relax\else Erratum \erratumnumber\nl\fi
Pages \startpage--\finishpage\nl}

\vglue 0.1truein   

{\parskip=0pt\leftskip 0pt plus 1fil\def\\{\par\smallskip}{\ifplaintex\large
\else\Large\fi\bf\thetitle}\par\medskip}   
\vglue 0.05truein 

{\parskip=0pt\leftskip 0pt plus 1fil\def\\{\par}{\sc\theauthors}
\par\medskip}%
 
\vglue 0.03truein 

{\small\leftskip 25pt\rightskip 25pt{\bf Abstract}\stdspace\theabstract

{\bf AMS Classification}\stdspace\theprimaryclass
\ifx\thesecondaryclass\relax\else; \thesecondaryclass\fi\par
{\bf Keywords}\stdspace \thekeywords\par}\vglue 7pt

}   

\font\phead=cmsl9 scaled 950
\font=cmsl9 scaled 1050
\font\pnum=cmbx10 scaled 913
\font\lnum=cmbx10 
\font\pfoot=cmsl9 scaled 950
\font=cmsl9 scaled 1050
\ifplaintex
\headline{\vbox to 0pt{\vskip -4.5mm\line{\small\phead\ifnum
\count0=\startpage ISSN 1464-8997 (on line)
1464-8989 (printed) \hfill {\pnum\folio}\else\ifodd\count0\def\\{ }%
\ifx\theshorttitle\relax\thetitle\else\theshorttitle\fi\hfill{\pnum\folio}
\else\def\\{ and }{\pnum\folio}\hfill\ifx\theshortauthors\relax\theauthors
\else\theshortauthors\fi\fi\fi}\vss}}
\footline{\vbox to 0pt{\vglue 0mm\line{\small\pfoot\ifnum\count0=\startpage
Published \publishdate:\qua\copyright\ \gtp\hfill\else
\gtm, Volume \thevolumenumber\ (\thevolumeyear)\hfill\fi}\vss
}}
\else
\makeatletter
\def\@oddhead{{\small\lhead\ifnum\count0=\startpage ISSN 1464-8997 (on line)
1464-8989 (printed) \hfill {\lnum\number\count0}\else\ifodd\count0
\def\\{ }\ifx\theshorttitle\relax \thetitle \else\theshorttitle\fi\hfill
{\lnum\number\count0}\else\def\\{ and }{\lnum\number\count0}
\hfill\ifx\theshortauthors\relax 
\theauthors\else\theshortauthors\fi\fi\fi}}\def\@evenhead{@oddhead}
\def\@oddfoot{\small\lfoot\ifnum\count0=\startpage Published \publishdate:\qua\copyright\ \gtp\hfill\else
\gtm, Volume \thevolumenumber\ (\thevolumeyear)\hfill\fi}
\def\@evenfoot{@oddfoot}
\makeatother
\fi

\let\maketitlepage\maketitlep
\let\makeshorttitle\maketitlepage
\let\maketitle\maketitlepage

\newwrite\gtoutfile
\long\gdef\makeheadfile{  
{\def\\{, }\def\s{ }
\immediate\openout\gtoutfile head.xxx
\immediate\write\gtoutfile{Proxy-for: \ifx\theasciiauthors\relax
\theauthors\else\theasciiauthors\fi\s<\ifx\theasciiemail\relax\theemail\else\theasciiemail\fi>}
\immediate\write\gtoutfile{\noexpand\\}
\immediate\write\gtoutfile{Authors: \ifx\theasciiauthors\relax
\theauthors\else\theasciiauthors\fi}
{\def\\{ }\immediate\write\gtoutfile{Title: \ifx\theasciititle\relax
\thetitle\else\theasciititle\fi}}
\immediate\write\gtoutfile{Subj-class: GT or SG, GR etc}
\immediate\write\gtoutfile{MSC-class: \theprimaryclass\ifx\thesecondaryclass\relax\else, \thesecondaryclass\fi}
\immediate\write\gtoutfile{Journal-ref: Geom. Topol. Monogr. \thevolumenumber\s
(\thevolumeyear) \startpage-\finishpage}
\immediate\write\gtoutfile{Comments: Published by Geometry and Topology Monographs at}
\immediate\write\gtoutfile{\s\s\s  http://www.maths.warwick.ac.uk/gt/GTMon\thevolumenumber/paper\thepapernumber.abs.html}
\immediate\write\gtoutfile{\noexpand\\}
\immediate\write\gtoutfile{}
\ifx\theasciiabstract\relax
\immediate\write\gtoutfile{\theabstract}\else
\immediate\write\gtoutfile{\theasciiabstract}\fi
\immediate\write\gtoutfile{}
\immediate\write\gtoutfile{\noexpand\\}
\immediate\write\gtoutfile{}
\immediate\closeout\gtoutfile}}  

\def\maketitlepage{\maketitlep\makeheadfile}
\let\makeshorttitle\maketitlepage
\let\maketitle\maketitlepage

\volumenumber{7}
\volumename{Proceedings of the Casson Fest}
\volumeyear{2004}
\papernumber{10}
\pagenumbers{235}{265}
\received{8 January 2004}
\revised{29 March 2004}
\accepted{16 March 2004}
\published{20 September 2004}

\usepackage{amsmath, amssymb, amscd}
\input rlepsf

\theoremstyle{plain}
\newtheorem{thm}{Theorem}[section]

\theoremstyle{definition}
\newtheorem{df}[thm]{Definition}

\def \P {\mathbb{P}}
\def \R {\mathbb{R}}
\def \Z {\mathbb{Z}}

\def \CA {\mathcal A}
\def \CB {\mathcal B}

\def \J {\mathcal J}

\def \CS {\mathcal S}
\def \CV {\mathcal V}
\def \CW {\mathcal W}

\def \a {\alpha}
\def \b {\beta}

\def \G {\Gamma}
\def \d {\delta}

\def \p {\phi}

\def \s {\sigma}
\def \Si {\Sigma}

\def \ti {\widetilde}

\def \bd {\partial}

\def \Q {Q}

\def \D {\Delta}

\begin{document}

\title[Ideal triangulations and spun normal surfaces]
{Ideal triangulations of 3--manifolds I:\\spun normal surface theory}

\authors{Ensil Kang\\J Hyam Rubinstein}
\address{Department of Mathematics, College of Natural 
Sciences\\Chosun University, Gwangju 501-759, Korea}
\secondaddress{Department of Mathematics and Statistics,
  The University of Melbourne\\Parkville, Victoria 3010, Australia}
\asciiaddress{Department of Mathematics, College of Natural 
Sciences\\Chosun University, Gwangju 501-759, Korea\\and\\Department 
of Mathematics and Statistics,
  The University of Melbourne\\Parkville, Victoria 3010, Australia}

\gtemail{\mailto{ekang@chosun.ac.kr}, 
\mailto{ekang@math.snu.ac.kr}{\rm\qua and\qua}\mailto{rubin@ms.unimelb.edu.au}}

\asciiemail{ekang@chosun.ac.kr, ekang@math.snu.ac.kr, rubin@ms.unimelb.edu.au}


\begin{abstract}
In this paper, we will compute the dimension of the space of spun
and ordinary normal surfaces in an ideal triangulation of the
interior of a compact 3--manifold with incompressible tori or Klein
bottle components. Spun normal surfaces have been described in
unpublished work of Thurston. We also define a boundary map from
spun normal surface theory to the homology classes of boundary
loops of the 3--manifold and prove the boundary map has image of
finite index. Spun normal surfaces give a natural way of
representing properly embedded and immersed essential surfaces in
a 3--manifold with tori and Klein bottle boundary  \cite{ka1},
\cite{ka2}.  It has been conjectured that every slope in a simple
knot complement can be represented by an immersed essential
surface \cite{ba}, \cite{bc}.  We finish by studying the boundary
map for the figure-8 knot space and for the Gieseking manifold,
using their natural simplest ideal triangulations. Some potential
applications of the boundary map to the study of boundary slopes
of immersed essential surfaces are discussed.
\end{abstract}

\asciiabstract{%
In this paper, we will compute the dimension of the space of spun and
ordinary normal surfaces in an ideal triangulation of the interior of
a compact 3-manifold with incompressible tori or Klein bottle
components. Spun normal surfaces have been described in unpublished
work of Thurston. We also define a boundary map from spun normal
surface theory to the homology classes of boundary loops of the
3-manifold and prove the boundary map has image of finite index. Spun
normal surfaces give a natural way of representing properly embedded
and immersed essential surfaces in a 3-manifold with tori and Klein
bottle boundary [E Kang, `Normal surfaces in knot complements' (PhD
thesis) and `Normal surfaces in non-compact 3-manifolds',
preprint].  It has been conjectured that every slope in a simple knot
complement can be represented by an immersed essential surface [M
Baker, Ann. Inst. Fourier (Grenoble) 46 (1996) 1443-1449 and (with D
Cooper) Top. Appl. 102 (2000) 239-252].  We finish by studying the
boundary map for the figure-8 knot space and for the Gieseking
manifold, using their natural simplest ideal triangulations. Some
potential applications of the boundary map to the study of boundary
slopes of immersed essential surfaces are discussed.}

\maketitle

{\small\it This paper is dedicated to Andrew Casson with thanks
  for his many wonderful mathematical ideas}\leftskip25pt\rightskip 25pt

\section{Introduction\label{Intro}}\leftskip0pt\rightskip 0pt
\label{sec1}

In a series of papers, starting with this one, we will investigate
ideal triangulations of the interiors of compact 3-manifolds with
incompressible tori or Klein bottle boundaries. Such
triangulations have been used with great effect, following the
pioneering work of Thurston \cite{th}. Ideal triangulations are
the basis of the computer program SNAPPEA of Weeks \cite{we} and
the program SNAP of Coulson, Goodman, Hodgson and Neumann
\cite{cg}. Casson has also written a program to find hyperbolic
structures on such 3-manifolds, by solving Thurston's hyperbolic
gluing equations for ideal triangulations.

In the second paper, we will study the problem of deforming a taut
triangulation \cite{la} to an angle structure. We show that given
a taut ideal triangulation of the interior of  a compact
3-manifold with tori or Klein bottle boundary components, which is
irreducible, $\P^2$-irreducible and atoroidal with no 2-sided
Klein bottles or Mobius bands,  there are angle structures if and
only if a certain combinatorial obstruction vanishes. This work is
inspired by ideas of Casson. In the third paper, various
connections between structures on ideal triangulations will be
examined in the case of small numbers of tetrahedra. In the fourth
paper, a theory of immersed normal and spun normal surfaces in
ideal triangulations with angle structures will be studied. In the
fifth paper, we show that immersed essential surfaces can be
homotoped to be spun normal, in case the surfaces do not lift to
fibres of bundles in finite sheeted covers. Examples are given of
different triangulations of bundles, where such fibres can and
cannot be realised by spun normal surfaces. Moreover, a very fast
algorithm to decide if a given embedded normal surface is
incompressible is given, using these ideas.

In this paper, for simplicity, all 3-manifolds $M$ will be the
interior of a compact manifold $N$ with tori or Klein bottle
boundary components.  We will work in the smooth category. All
3-manifolds will be irreducible and $\P^2$-irreducible, i.e every
embedded 2-sphere bounds a 3-ball and there are no embedded
2-sided projective planes. For basic 3-manifold theory, see either
\cite{he} or \cite{ja}.

An ideal triangulation $\G$ of $M$ will be a cell complex which is a
decomposition of $M$ into tetrahedra $\D_1, \D_2, ..., \D_k$ glued
along their faces and edges, so that the vertices of the tetrahedra
are all removed. Moreover the link of each such missing vertex will be
a Klein bottle or torus. Note that several vertices, edges or faces of
a tetrahedron may be identified. (Sometimes such triangulations are
called pseudo-triangulations.)  Using Moise's construction of
triangulations of 3--manifolds \cite{mo}, one can convert a
triangulation of $N$ into such an ideal triangulation, by collapsing
the boundary surfaces to ideal vertices and also collapsing edges
which join the ideal vertices to the interior vertices. See \cite{jr}
for a discussion of such collapsing procedures. One has to be careful
to ensure that at each stage of such collapsings the topological type
of $M$ is not changed.

Alternatively one can collapse $M$ onto a 2--dimensional spine
$\CS$ and then choose a dual triangulation. Assume that $M$ has at
least two boundary components. A convenient way to do this is to
choose a Heegaard splitting, ie a decomposition of $N$ into
compression bodies $X,Y$ so that these are glued along a Heegaard
surface $S$. $X,Y$ are obtained from $S \times I$ by gluing
2--handles to $S \times \{ 1 \}$ ($S$ can be non-orientable or
orientable). We can arrange that $S$ has boundary components of
$N$ on either side of it. Then the common boundary surface $S$ of
the compression bodies is obtained by gluing the two copies of $S
\times \{ 0\}$. Now if a collection of core disks for the
2--handles for $X,Y$ is attached to $S$, we get the required spine
$\CS$ for $M$ and $N$. It is straightforward to check that the
dual cell structure to this spine is an ideal triangulation for
$M$. If $M$ has only one boundary component, then the same
procedure will construct a triangulation of $M$ with one ideal
vertex and one vertex in the interior of $M$. We then (carefully)
collapse an edge joining these two vertices as in \cite{jr} to get
an ideal triangulation.

We now summarize Haken's theory of normal surfaces \cite{ha}, as
extended by Thurston to deal with spun normal surfaces in ideal
triangulations. Given an abstract tetrahedron $\D$ with vertices
$ABCD$, there are four normal triangular disk types, cutting off
small neighborhoods of each of the four vertices. There are also
three normal quadrilateral disk types, which separate pairs of
opposite edges, such as $AB, CD$. Each tetrahedron $\D_i$ of $\G$
contributes 7 coordinates which are the numbers $n_j$ of each of
the normal disk types. We can form a vector of length $7k$ from
all of these coordinates $n_j$, $ 1 \le j \le 7k$ and a normal
surface $S$ is formed by gluing finitely many normal disk types
together. There are $6k$ compatibility equations, each of the form
$n_i + n_j = n_m + n_p$, where the left side of the equation gives
the number of normal triangles and quadrilaterals with a
particular normal arc type in the boundary, eg the arc running
between edges $AB,AC$ in $\D$. If the face $ABC$ is glued to
$A'B'C'$ of the tetrahedron $\D'$, then $n_m, n_p$ are the number
of normal triangles and quadrilaterals with the boundary normal
arc type running between $A'B',A'C'$ in $\D'$. Note that we allow
self-identifications of tetrahedra and hence also of normal disk
types.

It turns out that the solution space $\CV$ of these compatibility
equations in $\R^{7k}$ has dimension $2k$, ie there are $k$
redundant compatibility equations. The non-negative integer
solutions in $\CV$ then give normal surfaces (possibly
non-uniquely) and we can regard $2k$ as the dimension of the space
of these surfaces. We will give a proof of this in the next
section, as it is important in computing the dimension of the
space $\CW$ of spun and ordinary normal surfaces. In fact, if $c$
is the number of tori and Klein bottle boundary components of $N$,
then we will show that the dimension of $\CW$ is $2k+c$. Note that
a normal surface is always closed, but may be immersed or have
branch points.

Before giving a definition of spun normal surfaces, we make some
preliminary remarks about ideal vertex linking surfaces. For
simplicity, in this
paper we are only dealing with the case of manifolds with either
boundary tori or Klein bottles
or both.
However much of what is done can be carried over to the case of
higher-genus boundary surfaces
and we will return to this topic in a later paper.

If $M$ is the interior of a compact manifold with boundary tori
and Klein bottles, with an ideal triangulation $\G$, then we will
denote by $T$ an ideal vertex linking torus or the orientable
double cover of an ideal vertex linking Klein bottle. For a given
homology class of a simple closed curve on $T$, we will choose a
representative closed curve $C_T$ with minimal length in the
1--skeleton of $T$, which is viewed as a normal surface
constructed from triangular disk types in $\G$. There is a
uniquely defined infinite cyclic covering space of $T$
corresponding to the subgroup of $\pi_1(T)$ generated by the
homotopy class of $C_T$. Let $\widetilde{C}$ denote a lift of $C_T$ to
this covering space. Although $C_T$ may not be simple, since it
has been chosen to have minimal length, it is easy to see that
$\widetilde{C}$ is a simple closed curve. So we find two half-open
infinite annuli bounded by $\widetilde{C}$, which are each covered by
lifted triangular disk types from $T$. Let $A_+, A_-$ denote the
images of these annuli, viewed as half-open infinite spiralling
annuli in $M$. Note these annuli have boundary the curve $C_T$.
For convenience, we label the tori boundary components and the
orientable double covers of the Klein bottle boundary components
of $M$ by $T_j$, for $1 \le j \le c$, and choices of simple closed
curves in $T_j$, for $j$ fixed, by ${C^i}_{T_j}$, where $ 1 \le i
\le m_j$.

Our definition of spun normal surfaces follows the spirit of
Thurston, emphasizing the geometric picture in the 3--manifold $M$.
However for calculations and applications, it turns out to
be more convenient to use $\Q$--theory as discussed below.

\begin{figure}[ht!]
\centering
\epsfxsize=4.5in
\epsfbox{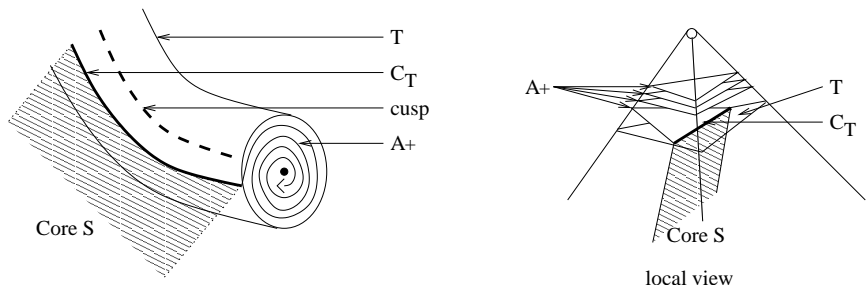}
\caption{Spun normal surface}
\label{fig1}
\end{figure}

\begin{df}
Suppose that $M$ is the interior of a compact manifold with tori
or Klein bottle boundary components and $\G$ is an ideal triangulation
of $M$. A spun normal surface $S$ has infinitely many normal disk
types, of which
finitely many are quadrilaterals and infinitely many are triangular.
Moreover $S$ has a finite collection of boundary slopes, which are choices
of closed curves ${C^i}_{T_j}$ as above, in the 1--skeleton of each
boundary torus component $T_j$ of $M$, or in the orientable
double cover $T_j$ of any Klein bottle component. We can glue
together all except a finite number of the triangular
disk types in $S$ to form a finite collection of half open spiralling
annuli ${A^{ij}}_+$ or ${A^{ij}}_-$ as described above.
Then the remaining finitely many quadrilateral and triangular disk
types satisfy the compatibility equations, together
with this finite collection of annuli which are viewed as having
boundary curves given by the normal arc types
of the curves ${C^i}_{T_j}$.
\end{df}

\begin{figure}[ht!]
\centering
\epsfxsize=4.5in
\epsfbox{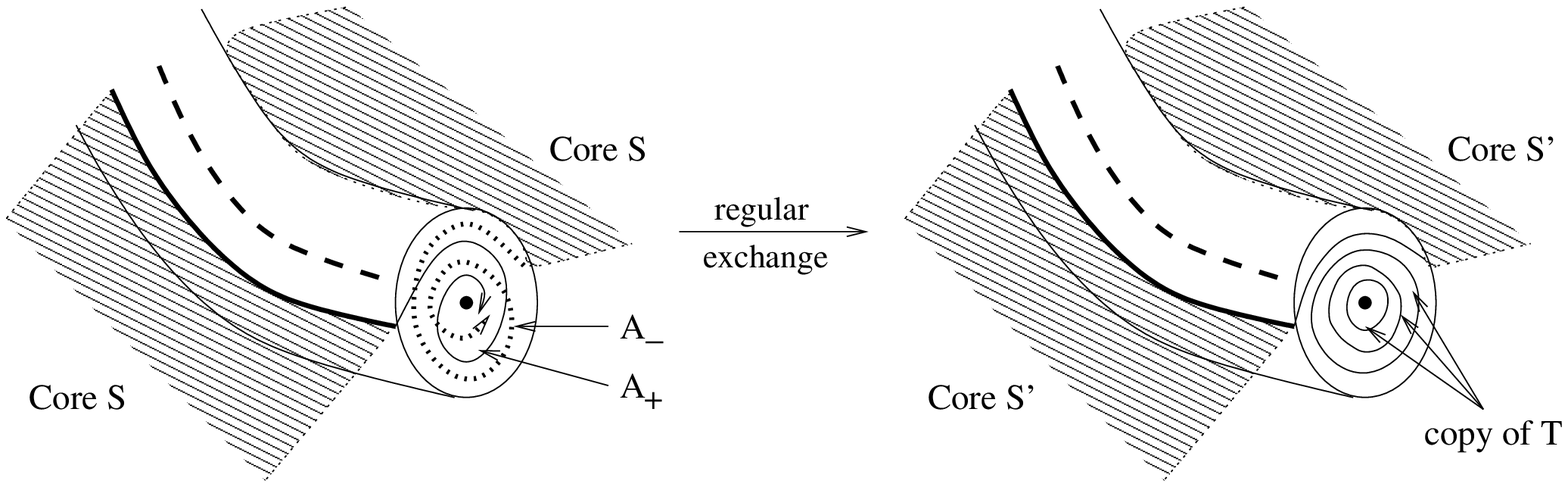}
\caption{Extra or cancelling infinite half open annuli}
\label{fig2}
\end{figure}

Notice that we can have `extra' or cancelling infinite half-open
annuli at the same boundary component of $M$, since if a spun normal
surface
$S$ has copies of both $A_+$ and $A_-$ for a given boundary slope,
then we can find a new spun normal surface
or possibly just a normal surface $S'$, by removing (or cancelling)
these two annuli from $S$. Moreover $S$ might
be embedded, immersed or have branch points. In the fifth paper of
our sequence, it will be natural to
consider also higher genus boundary surfaces for $M$. In that case, we must
consider generally immersed spiralling annuli $A_+,A_-$, as occur in
the case of Klein bottle boundary components.
These arise naturally, by projecting embedded annuli in a suitable
covering space of a boundary surface $S$,
to $M$.

It will be very convenient to extend our ideas of normal and spun
normal surfaces, to allow negative numbers of
triangular and quadrilateral disk types. For normal surfaces, this is
easy; we just take all integer solutions of the compatibility
equations and call these formal normal surfaces,  to allow
coordinates (multiples of the disk types) to be any integers, not
just non-negative ones.
For spun normal surfaces we can play the same game by allowing the
half open annuli to occur with negative multiplicity as well
as the disk types and have formal solutions of the compatibility
equations as well. We will refer to these as formal spun normal
surfaces.

Now to form a vector space $\CW$ of formal spun and ordinary normal
surfaces $S$, we will consider only the quadrilateral coordinates
of each $S$. So $\CW$ will be a subspace of $\R^{3k}$. This idea
has been studied previously in \cite{to}, in case of ordinary
normal surface theory in standard triangulations and is called
$\Q$--normal surface theory. For spun normal surfaces, $\Q$--theory
has been investigated in \cite{ka1},\cite{ka2}. There are $k$
compatibility equations for the quadrilaterals and we discuss this
in detail later on. In an ideal triangulation, the solutions to
these equations are naturally either normal or spun normal
surfaces. The only surfaces which are not `seen' by this theory,
are multiples of the boundary Klein bottles and tori, formed
entirely of triangular disk types. If we added these in also, the
theory would have dimension $2k+2c$. However these boundary
surfaces play no important role, so it is reasonable to leave
them out of consideration. Spun normal surface theory has been
used in an interesting way by Stefan Tillmann \cite{ti}, to study
essential splitting surfaces arising from representation varieties
in Culler--Shalen theory.

In our fifth paper, we will prove an extension of an unpublished
result of Thurston,
namely that any immersed essential surface is homotopic to a spun
normal surface, unless
the surface lifts to a fibre of a bundle structure in some finite
sheeted covering space.
Consequently, if a certain boundary slope is not in the image of the
boundary map
(as in the examples in section 5 below) then either there are no
immersed essential surfaces
with such a boundary slope (in the homological sense of our boundary
map), or any such
a surface must be a virtual fibre.

Note that for the examples of the figure--8 knot complement and the
Gieseking manifold in this paper, we are not able to get any
interesting information on the boundary slopes of their immersed
incompressible surfaces. These two examples are instances of
layered triangulations of bundles, and the latter give special
problems for the question of homotoping immersed incompressible
surfaces to be spun normal. Full details are in the fifth paper.

\section{The canonical basis for normal surface theory}
\label{sec2}

In this section, we will construct a very useful `canonical basis'
for the solution space $\CV$ of the compatibility equations for
formal normal
surface theory. This is independent of the details of the
triangulation and the 3--manifold and was inspired by ideas of Casson, although
he never stated a result like this. The theorem plays a key role in
the second paper in our
sequence, where we study the problem of deforming taut structures on
ideal triangulations to angle structures.

\begin{thm}

Let $M$ be a 3--manifold with a triangulation $\G$ which is either
ideal or not, depending on whether $M$ is the interior of a compact
manifold $N$ with tori and Klein bottle boundary components, or $M$ is
closed. Assume there are
$k$ tetrahedra and $e$ edges in $\G$. Then
$\CV$ has a basis consisting of $k$ formal `tetrahedral' surfaces and $e$
formal `edge' surfaces.
\end{thm}

\begin{proof}
The idea is to study $e$ linear functionals
$\p_1, \p_2, ... , \p_e$ from $\R^{7k}$ to $\R$, one for each edge
$E_1, E_2, ... , E_e$ of $\G$. Each of these functionals restricts
to $\CV$ and we obtain $e$ linearly-independent functionals, which
will be denoted by the same symbols. It then follows that $\CV$
splits as a direct sum $\CA \oplus \CB$, where $\CA$ is the common
kernel of all the linear functionals and $\CB$ is a `natural'
cokernel. We will then find bases of $\CA, \CB$ separately and
these will form our canonical basis.

For an edge $E_i$, let $d_i$ denote its degree, ie the number of
copies of edges of tetrahedra which are glued together to form
$E_i$. For every normal triangle or quadrilateral with $n$
vertices on $E_i$, we map this normal disk type to $\frac{n}{d_i}$.
Extending by linearity gives a linear functional $\p_i$ from
$\R^{7k}$ to $\R$. To show these linear functionals restricted to
$\CV$ are linearly independent, we construct a collection of
vectors $\b_1, \b_2, ..., \b_e$ in $\CV$, so that $\p_i(\b_j) =
-2\d_{ij}$, where as usual $\d_{ij}=0$ if $i \ne j$ and
$\d_{ii}=1$. This will also show that $\b_1, \b_2, ..., \b_e$ is a
basis for $\CB$ as required. We call these vectors `edge'
solutions.

Consider the collection of triangular disk types which have at least
one vertex on an edge $E_i$. We take $-n$ copies of each
such a triangular disk, where the disk has precisely $n$ vertices on
$E_i$, for $1 \le n \le 3$. Similarly consider the set of
tetrahedra which contain at least one edge which becomes $E_i$ after
the faces and edges are glued together.
Suppose $ABCD$ represents such a tetrahedron before gluing and $AB$
is an edge which becomes $E_i$. Then one copy of the quadrilateral
separating
$AB$ from $CD$ is selected. Notice that if the edge $CD$ also becomes
$E_i$, then we must choose two copies of this quadrilateral.
Do this for the $d_i$ edges of tetrahedra which are glued together to
form $E_i$. Then $\b_i$ is defined as this collection of negative
triangular
disks and positive quadrilateral disks.

We need to check that $\b_i$ is indeed in $\CV$, ie is a solution
of the compatibility equations. This is straightforward -- consider
a face $ABC$ of a tetrahedron $ABCD$, where $AB$ is identified
with $E_i$. For the normal arc type running between $AC$ and $AD$,
there is a triangular disk (cutting off $A$) and a quadrilateral
disk (separating $AB$ and $CD$) which are `glued' along this arc.
Notice the triangular disk is taken with sign $-1$ whereas the
quadrilateral is taken with a positive sign. Hence the boundary
arc types cancel out as required by the compatibility equations.
Assume next that the tetrahedron $A'B'C'D'$ is glued to $ABCD$ by
identifying $ABC$ with $A'B'C'$. For the normal arc type running
between $AB$ and $AC$, the triangular disk type cutting off $A'$
is glued to the triangular disk cutting off $A$, along this arc.
Since both disks have sign $-1$, the boundary arcs cancel.
Finally, a similar argument applies for the normal arc type
running between $AC$ and $BC$ by gluing two quadrilateral disk
types, separating $AB$ and $CD$ (respectively $A'B'$ and $C'D'$)
in the two tetrahedra. One can check that self-identification of
tetrahedra does not affect this cancellation and so $\b_i$ is a
`formal' normal surface, ie a solution of the compatibility
equations. If there are no tetrahedra with two or more edges
identified to $E_i$, the surface $\b_i$ can be visualized as a
cylindrical box with negative top and bottom and positive sides.

Next we observe that $\p_i(\b_j) = -2\d_{ij}$. It is clear that
$\p_i(\b_i) = -2$. If there are no tetrahedra
with two or more edges identified to $E_i$, then the functional
$\p_i$ takes values $-1$ at the top and bottom ends of the box and $0$
at all other vertices. So the total value of $\p_i(\b_i)$ is $-2$ as claimed.
Again one can check that self-identifications do not affect this
calculation. It remains to show that $\p_j(\b_i) = 0$ if $i \ne j$.
(We have interchanged $i,j$ to use our previous description of the
solution $\b_i$).
With notation as in the previous paragraph, assume that the edge
$E_j$ is $AC$. Now there are two triangular disks with signs $-1$ and
two
quadrilateral disks with signs $+1$ coming from $\b_i$, at this edge, so
their contributions to $\p_j(\b_i)$ cancel out. Self-identifications will
not alter this
phenomenon and so we have proved both that the surfaces $\b_1, \b_2,
..., \b_e$ are linearly independent and also that they form a basis
for
a cokernel $\CB$ for the linear functionals $\p_1, \p_2, ... , \p_e$
on $\CV$.

To complete the proof of this theorem, we need to construct a basis
$\a_1, \a_2, ...,\a_t$ for the kernel $\CA$ of the linear functionals.
Again let $ABCD$ denote the $i$th tetrahedron. Then $\a_i$ will be
each of the $4$ triangular disk types cutting off a vertex taken with
sign
$-1$ and each of the $3$ quadrilateral disk types with sign $+1$ in
this tetrahedron. We must check that $\a_i$ is in $\CV$, ie
solves the compatibility equations. This is very similar to the
previous argument - for a normal arc type between edges $AB$ and
$AC$, the
triangular disk cutting off $A$ has sign $-1$ and the quadrilateral
disk separating $AB$ and $CD$ has sign $+1$, so the two boundary arcs
coming from these disks cancel as required.

It is equally easy to see that all the surfaces $\a_1, \a_2,
...,\a_t$ are in the common kernel $\CA$ of the linear functionals
and that they are linearly independent, since they have no
non-zero coordinates in common. The final step is to prove that
these surfaces span $\CA$. So we need to consider any vector $S$
which is in $\CA$ and show it is a linear combination of $\a_1,
\a_2, ...,\a_t$. Notice that the `intersection number' of $S$ with
every edge is zero, since $S$ is in $\CA$. If as usual, $ABCD$ is
the $i$th tetrahedron of $\G$, then the compatibility equations
into a neighboring tetrahedron $A'B'C'D'$, where the face $ABC$ is
glued to $A'B'C'$, show that the total number $N$ (which can be
positive or negative) of triangle and quadrilateral disks in
$ABCD$ which meet $AB$ is the same as the total number of triangle
and quadrilateral disks in $A'B'C'D'$ which meet $A'B'$. Following
the tetrahedra around the edge $E_i$, which is $AB$ and $A'B'$
after identifications, we see this number $N$ is independent of
the tetrahedron containing an edge identified to $E_i$. Since the
value of the functional $\p_i$ is exactly $Nd_i$, where $d_i$ is
the degree of $E_i$, it follows that $N = 0$, since $S$ is in
$\CA$. It is an easy exercise to check that the only sets of
numbers of triangles and quadrilaterals in $ABCD$ which satisfy
this condition for each of the $6$ edges of the tetrahedron, are
multiples of $\a_i$. Hence this proves that $\a_1, \a_2, ...,\a_t$
span $\CA$ and the construction of the canonical basis is
complete.
\end{proof}

{\bf Remarks}\qua
For ideal triangulations, where the 3--manifold $M$ is the interior
of a compact 3--manifold with tori and Klein bottle boundaries, a
simple Euler characteristic argument, using Poincar\'e duality,
shows that $k = e$, ie the number of tetrahedra is equal to the
number of edges. This is the case of most interest to us. Another
important case is a triangulation of a closed 3--manifold with a
single vertex \cite{jr}. Again a simple Euler characteristic
argument shows then that $e = k+1$. So we conclude that the
dimension of $\CV$ is either $2k$ or $2k+1$ respectively in these
two standard cases.

\section{Dimension of spun and ordinary normal surface theory}
\label{sec3}

In this section, we will give a lower bound on the dimension of
formal spun normal surface theory, working with an ideal
triangulation of $M$. We begin with a quick summary of $\Q$--normal
surface theory (\cite{to}, \cite{ka1}, \cite{ka2}). Let $m_1, m_2,
...,m_{3k}$ denote the number of quadrilateral normal disk types
in an ordinary or spun normal surface. It is not hard to show that
once the quadrilateral disks are specified, a formal normal or spun
normal surface can be reconstituted, up to multiples of the
boundary Klein bottles and tori. The key fact we need is that the
quadrilaterals satisfy $k = e$ compatibility equations, one for
each edge of $\G$. These equations are readily found by
eliminating the triangular coordinates from the ordinary
compatibility equations in the previous section. Before
introducing these equations, we note an important way of labelling
the corners of the quadrilateral disk types, assuming that $M$ is
orientable.

\begin{figure}[ht!]
\cl{
\relabelbox\small
\epsfxsize=3.2in
\epsfbox{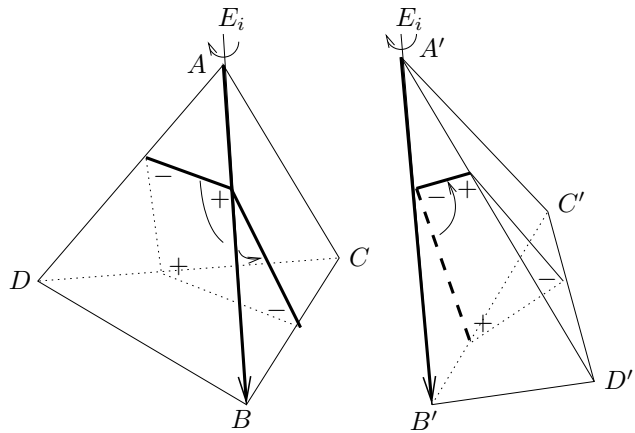}
\relabel{A}{$A$}
\adjustrelabel <+5pt, 0pt> {B}{$B$}
\relabel{C}{$C$}
\adjustrelabel <-3pt, 0pt> {D}{$D$}
\relabel{A'}{$A'$}
\adjustrelabel <-15pt, 2pt> {B'}{$B'$}
\relabel{C'}{$C'$}
\relabel{D'}{$D'$}
\relabel{E1}{$E_i$}
\relabel{E2}{$E_i$}
\adjustrelabel <-3pt,-1pt> {+1}{$+$}
\adjustrelabel <-1pt,0pt>  {+2}{$+$}
\adjustrelabel <0pt,+2pt>  {+3}{$+$}
\adjustrelabel <-3pt,0pt>  {+4}{$+$}
\adjustrelabel <-1pt,-1pt> {-1}{$-$}
\adjustrelabel <-2pt,+1pt> {-2}{$-$}
\relabel{-3}{$-$}
\adjustrelabel <-2pt, 1pt> {-4}{$-$}
\endrelabelbox}
\caption{Labelling the corners of quadrilateral disk types}
\label{fig3}
\end{figure}

Let $ABCD$ and $A'B'C'D'$ be tetrahedra which are glued along the
faces $ABC$ and $A'B'C'$. Let $AB$ and $A'B'$ be identified with the
$i$th edge $E_i$ of $\G$. For the quadrilateral separating $AC$ and
$BD$, we will associate corner signs, $+$ to the corners on $AB, CD$
and $-$ to corners on $AD, BC$ (see Figure \ref{fig3}).  Similarly for
the quadrilateral separating $AD$ and $BC$, $+$ is given to the
corners on $AC, BD$ and $-$ to corners on $AB, CD$.  In the adjacent
tetrahedron $A'B'C'D'$, the signs are reversed. So for the
quadrilateral separating $A'C'$ and $B'D'$ (respectively $A'D'$ and
$B'C'$), we will associate corner signs, $-$ to the corners on $A'B',
C'D'$ (respectively $A'C', B'D'$) and $+$ to corners on $A'D', B'C'$
(respectively $A'C', B'D'$) (see Figure 3). To show that this
assignment of signs is consistent for all quadrilaterals
simultaneously, we note that it is equivalent to introducing an
orientation on each tetrahedron. For example, a right handed rule on
the tetrahedron $ABCD$ can be viewed as having orientations from $A$
to $B$ on the edge $AB$ and then circulating from the face $ABD$ to
the face $ABC$ so that the normal arc running between $AB$ and $AD$
cutting off $A$ is transformed into the normal arc running between
$AB$ and $BC$ cutting off $B$ in the quadrilateral separating $AC$ and
$BD$ (see Figure \ref{fig3}).

\noindent We associate this with the positive sign in the corner
of the quadrilateral at the corner on $AB$. It is easy to see that
if the orientation of $AB$ is reversed, we still have a positive
sign at this corner coming from the right hand rule. By
translating a right handed system throughout $M$, we see that the
consistency of signs is guaranteed. Now for an orientable
manifold, we can define a quadrilateral compatibility equation by
taking the sum of all quadrilateral types meeting a fixed edge
$E_i$ with sign given by the corner of each quadrilateral at
$E_i$. So half the quadrilaterals have sign $+$ and half have sign
$-$.

On the other hand if $M$ is non-orientable, then choosing an
orientation reversing loop passing through centers of tetrahedra
and faces, following around the quadrilaterals in these tetrahedra
will give an inconsistency in this choice of corner signs. However
we can still give a sign to corners of quadrilaterals around each
fixed edge $E_i$, but cannot do this consistently for all edges
and corners simultaneously. We do obtain a similar compatibility
equation in the non-orientable case, since it does not matter how
the quadrilateral corner signs vary from one edge to the other.

We are now ready to give the main result of this section.

\begin{thm}
Let $M$ be the interior of a compact manifold with boundary components
which are tori or Klein bottles. Let $\G$ denote an ideal triangulation of
$M$ and let $\CW$ be the vector subspace of formal normal and spun
normal surfaces. Then
the dimension of $\CW$ is at least $2k+c$, where $k$ is the number of
edges in $\G$ and $c$ is the number of boundary
components.

\end{thm}

\begin{proof}
We first show that there is at least one redundancy in the set
of $\Q$--matching equations, if $M$ is orientable. For adding all
the equations, each quadrilateral appears in at most 4 equations
with a total of exactly two $+$ signs and two $-$ signs. Consequently the
entries corresponding to this quadrilateral will cancel, ie the sum of
all $\Q$--matching equations is zero. This proves that we need at
most $k-1$ matching equations and the dimension of $\CW$ is at
least $2k+1$.

\begin{figure}[ht!]
\centering
\relabelbox\small
\epsfxsize=5.1in
\epsfbox{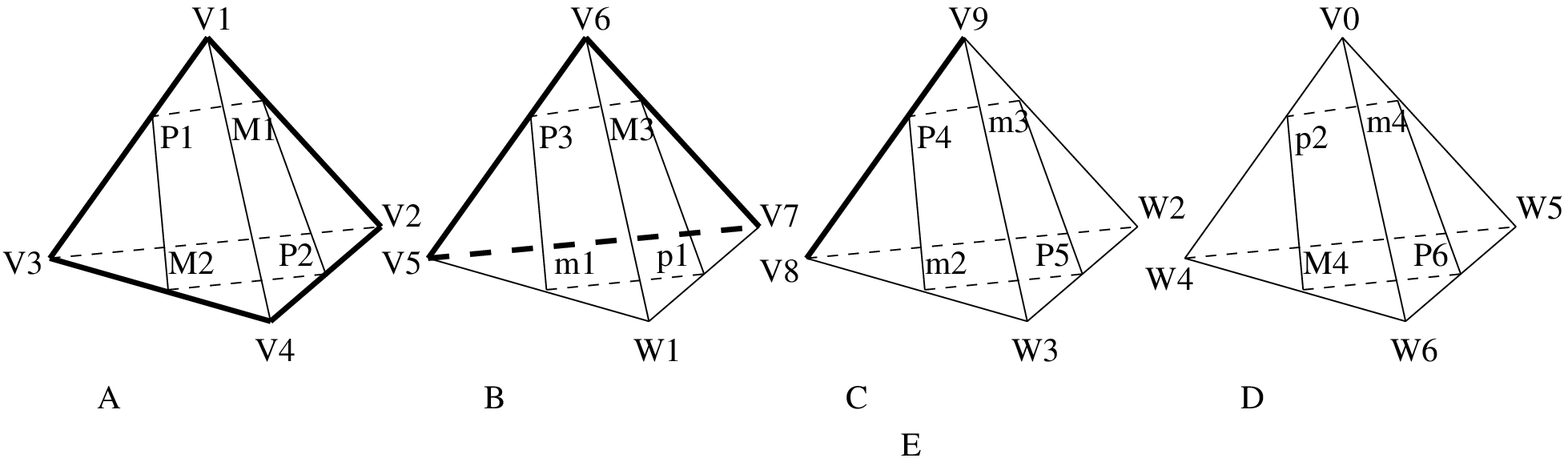}
\relabel{V1}{$V_1$}
\relabel{V2}{$V_1$}
\relabel{V3}{$V_1$}
\relabel{V4}{$V_1$}
\relabel{V5}{$V_1$}
\relabel{V6}{$V_1$}
\relabel{V7}{$V_1$}
\relabel{V8}{$V_1$}
\relabel{V9}{$V_1$}
\relabel{V0}{$V_1$}
\relabel{W1}{$V_2$}
\relabel{W2}{$V_2$}
\relabel{W3}{$V_2$}
\relabel{W4}{$V_2$}
\relabel{W5}{$V_2$}
\relabel{W6}{$V_2$}
\relabel{P1}{$\oplus$}
\relabel{P2}{$\oplus$}
\relabel{P3}{$\oplus$}
\relabel{P4}{$\oplus$}
\relabel{P5}{$\oplus$}
\relabel{P6}{$\oplus$}
\relabel{p1}{$+$}
\relabel{p2}{$+$}
\relabel{M1}{$\ominus$}
\relabel{M2}{$\ominus$}
\relabel{M3}{$\ominus$}
\relabel{M4}{$\ominus$}
\relabel{m1}{$-$}
\relabel{m2}{$-$}
\relabel{m3}{$-$}
\relabel{m4}{$-$}
\adjustrelabel <0pt, -5pt> {A}{$\sum\varepsilon(V_1V_1){=}0$}
\adjustrelabel <-5pt, -5pt> {B}{$\sum\varepsilon(V_1V_1){=}0$}
\adjustrelabel <-15pt,-5pt> {C}{$\sum\varepsilon(V_1V_1){=}\sum\varepsilon(V_2V_2)$}
\adjustrelabel <10pt, -5pt> {D}{$\sum\varepsilon(V_2V_2){=}0$}
\endrelabelbox
\caption{The cases of two ideal vertices}
\label{fig4}
\end{figure}

The next step is to analyze the case when there is more than one
boundary torus, supposing that $M$ is orientable. So there are
several ideal vertices, $V_1, V_2, ...,V_c$. Now label any
tetrahedron $ABCD$ by specifying which ideal vertices are
identified with $A,B,C,D$. Let us consider first the case that
there are only two ideal vertices $V_1, V_2$ (see Figure \ref{fig4}). There
are 5 cases, depending on how many of the vertices $A,B,C,D$ are
labelled $V_1$. We prove that adding all the matching equations
corresponding to edges with both ends labelled $V_1$ gives the
same result as adding over all edges with both vertices at $V_2$.
So there is a second redundancy in the matching equations and the
dimension of $\CW$ is at least $2k+2$. Note that if the number $n
\ne 2$ of vertices in some tetrahedron which are labelled $V_1$,
then adding the appearances of a quadrilateral type in this
tetrahedron over all matching equations corresponding to edges
with both ends at $V_1$ will give zero. This is because the signs
at the corners of the quadrilateral at such edges always cancel.
For the final case of two vertices labelled $V_1$ and two labelled
$V_2$, we see that the number of appearances of a quadrilateral
type in the two sums of matching equations is the same.
Consequently this case is complete.

Generally, if there are ideal vertices, $V_1, V_2, ...,V_c$, then we
find (by the same argument) that the sum of the matching equations at
edges labelled $V_1V_1$ is equal to the sum of the matching equations
at all edges labelled $V_iV_j$ for $i,j>1$ and similarly for $V_2V_2,
...,V_cV_c$. Note that by connectivity of the one-skeleton of $\G$, we
have paths of edges between any two ideal vertices. This is easily
seen to imply that there are at least $c-1$ such redundant matching
equations and one redundancy that the sum of all $\Q$--matching
equations is zero. Hence the dimension of $\CW$ is at least $2k+c$. In
fact, one can consider a collection of variables labelled $x_{ij}$,
one for each equation which is a sum of the compatibility or matching
equations at all edges labelled $V_iV_j$. Then our relations found
above can be written as $x_{mm} = \Si \{x_{ij}:i,j \ne m\}$ for
$m=1,2,\ldots,c-1$ and $\Si x_{ij} = 0$. Now it is easy to put this
system of linear equations into row-echelon form and check that the
rank is indeed $c$, proving that there are at least $c$ redundant
compatibility equations.

Assume next that $M$ is non-orientable. Let $\ti M$ denote the
orientable double cover of $M$ with $\ti \G$ the lifted
triangulation. There are pairs of edges interchanged by the
covering involution $\s$. We are interested in vectors of
quadrilaterals which are invariant under the induced action of
$\s$, since any lifted normal or spun normal surface will lift to
such a vector. So we introduce $3k$ new independent variables, by
forming a copy of $\R^{3k}$ embedded in $\R^{6k}$, given by $\a
\rightarrow (\a, \s (\a))$ for a vector $\a$ in $\R^{3k}$. Our
matching equations can be thought of as being located around two
edges $E, \s(E)$ simultaneously. This is achieved by using the new
variables, which are integer multiples of pairs of quadrilateral
disks $(D, \s(D))$. As $\s$ reverses orientation, it must reverse
our choice of signs of corners of quadrilateral disks in $\ti M$.
Hence corresponding quadrilateral coordinates in $\ti M$ under the
induced action of $\s$ are negatives of each other, ie $n(D,
\s(D))$ corresponds to the pair of coordinates $(n,-n)$ in $\ti
M$.

The conclusion is that we do get a similar theory to the orientable case,
starting with $k$ matching equations and finding $c$ redundancies,
proving the dimension estimate is also valid in the non-orientable case.
Note that $c$ is the number of Klein bottle and tori boundary components
of $M$, since each Klein bottle lifts to a torus which is invariant
under $\s$ and a torus lifts to a pair of tori interchanged by $\s$. So
since invariant boundary components or pairs of switched components
give rise to redundancies, this is the correct definition of $c$.

To illustrate this, consider the simple case that $M$ has a single
boundary component which is a Klein bottle or torus. This lifts to a
$\s$--invariant torus or a pair of tori interchanged by $\s$ in $\ti M$.
We can label the ideal vertex or vertices by this torus or tori. The
single relation amongst our new pairs of quadrilateral coordinates
becomes just that the sum of all the matching equations is zero. If there
were more boundary components, then an analogous argument to the previous
one shows that there is at least one redundancy amongst the matching
equations of the pairs of quadrilateral coordinates for each boundary
component of $M$.

Next, if there are two boundary components of $M$, then we get two
classes of ideal vertices of $\ti M$. Each class is either a single
vertex which is fixed under the action of $\s$ or a pair of vertices
interchanged by $\s$. As previously, it is easy to check with our new
coordinates, that the sum of all the matching equations for edges with
both ends in one vertex class gives the same result as the sum over all
edges with both ends in the other class. The general case follows the
same pattern as before.
\end{proof}

\section{The boundary map for spun normal surface theory}
\label{sec4}

To complete the computation of the dimension of $\CW$, we need to
prove it is also bounded above by $2k+c$. To do this,
we will define a boundary map $\bd\co \CW \rightarrow
\R^{2c}$, where the range of $\bd$ is the direct sum of all the
homology groups $H_1(T_i,\R)$, where the $T_i$
are the tori boundary components or the orientable double covers
of the Klein bottle boundary components of $M$, for $1 \le i \le c$.
  A key property of $\bd$ is that its kernel $A$ is
precisely the quadrilateral vectors corresponding to all formal normal
surfaces in $\CW$. Hence we can compute readily the dimension of $A$ as
$2k-c$, since the natural projection
$\CV \rightarrow \CW$, which erases the triangular coordinates of a
normal class, has kernel generated by the boundary tori and Klein
bottles, so has dimension $c$. Since the dimension of $\CV$ is
$2k$, we will see that the subspace $A$ of standard normal classes in
$\CW$ has dimension $2k-c$. Consequently, the dimension of $\CW$
is at most $2k-c +2c=2k+c$ and the computation of dimension will be
complete. Moreover we obtain the important result that $\bd$ has
image of dimension $2c$, so maps onto a subgroup of finite index, when
using $\Z$ coordinates. We compute the image of $\bd$ for the figure--8
knot space and Gieseking manifolds in the final section.

\begin{thm}
Let $M$ be the interior of a compact manifold with boundary components
which are tori or Klein bottles. Let $\G$ denote an ideal triangulation of
$M$ and let $\CW$ be the vector subspace of formal normal and spun
normal surfaces. Then
the dimension of $\CW$ is $2k+c$, where $k$ is the number of edges in
$\G$ and $c$ is the number of boundary
components. Moreover the boundary map $\bd\co \CW \rightarrow
\R^{2c}$ is onto and hence when using $\Z$ coefficients, it follows
that $\bd$ has image of finite index
in the direct sum $\Z^{2c}$ of the homology groups $H_1(T_i,\Z)$,
where the $T_i$
are the tori boundary components or the orientable double covers
of the Klein bottle boundary components of $M$, for $1 \le i \le c$.
\end{thm}

\begin{proof}
As in the previous section, we begin with the case that $M$ is
orientable, so we can put signs on all corners of quadrilateral
disks in a consistent manner. We will then proceed to put
orientations on all boundary arcs of quadrilaterals. As
previously, let $ABCD$ denote a tetrahedron and consider the
quadrilateral separating edges $AC$ and $BD$. Assume that a $+$
sign is attributed to the vertex on the edge $AB$ as in the
previous section. We now orient the boundary arc of this
quadrilateral running between the vertices $u, u'$ on $AB, AD$
respectively with an arrow from $u$ to $u'$ (see Figure \ref{fig5}). If
$v,v'$ are the vertices on the edges $BC,CD$ respectively, then
the boundary edges are oriented in the directions $v'v, uv, v'u'$,
ie from a vertex of $+$ sign to a vertex of $-$ sign. Notice that
if we adopt the convention that a triangular face has the
tetrahedron behind it, with a normal boundary arc having the
vertex cut off above it and the positive corner is on the left and
the negative on the right, then the orientation of the arc is from
left to right. All $4$ arcs are now oriented the same way by this
convention, as are all arcs of all quadrilaterals.

\begin{figure}[ht!]
\cl{\relabelbox\small
\epsfxsize=2.2in
\epsfbox{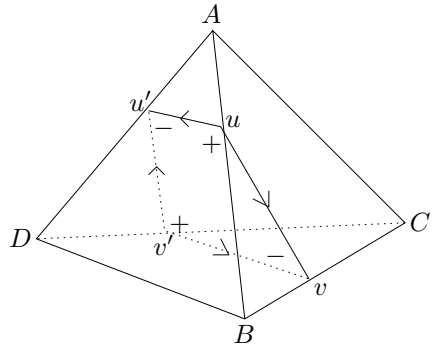}
\adjustrelabel <-2pt,0pt> {A}{$A$}
\adjustrelabel <-2pt,0pt> {B}{$B$}
\adjustrelabel <-2pt,0pt> {C}{$C$}
\adjustrelabel <-2pt,0pt> {D}{$D$}
\adjustrelabel <-2pt,0pt> {u}{$u$}
\adjustrelabel <-2pt,0pt> {u'}{$u'$}
\adjustrelabel <-2pt,0pt> {v}{$v$}
\adjustrelabel <-2pt,0pt> {v'}{$v'$}
\adjustrelabel <-2pt,0pt> {+}{$+$}
\adjustrelabel <-2pt,0pt> {-}{$-$}
\adjustrelabel <-2pt,0pt> {+1}{$+$}
\adjustrelabel <-2pt,0pt> {-1}{$-$}
\endrelabelbox}
\caption{Orientation of the boundary arcs of a quadrilateral}
\label{fig5}
\end{figure}

Next we want to glue together quadrilaterals at corner vertices
and possibly along some boundary arcs. Imagine two `adjacent' quadrilaterals
with such a vertex in common at some edge of $\G$, eg the corner
$u$ at the edge $AB$ in the previous paragraph for the first
quadrilateral, chosen so that there is a wedge of triangular disks
at $u$ between these two quadrilaterals on a (possibly spun)
normal surface $S$. There could be an empty wedge, if these two
quadrilaterals are actually glued along an edge, eg the edge $uu'$,
rather than just
glued at the corner vertex $u$. Now it
is easy to see the sign of the corner of the second quadrilateral
at $u$ must be $-$, ie the signs of the two adjacent quadrilaterals at
such a vertex are always opposite. So we proceed to glue together
quadrilaterals with vertices on the same edge of $\G$ with
opposite signs at these vertices.

Notice by the compatibility equations for quadrilaterals,
eventually at $u$ we find an even number of quadrilaterals with
alternating signs glued together in this fashion. Do this for all
vertices of all quadrilaterals in any formal vector of
quadrilaterals in $\CW$. We call this the quadrilateral part of
the surface and denote it by $S_Q$. There are finitely many ways
of producing a quadrilateral part of a normal surface, which can
then be filled in by triangular disks to form $S$. If the boundary
curves of $S_Q$ are inessential on the boundary tori of $M$, then
we can complete $S$ using disks and once-punctured tori, and
obtain a closed normal surface. If these curves are essential, we
get a spun normal surface. See \cite{ka1}, \cite{ka2} for further
details on this.

Now the key point is that since our choice of orientations of
boundary arcs is independent of the gluings, and when two
quadrilaterals are glued
along such an arc, then the two orientations cancel, we see that the
total homology class of the boundary curves of $S$ is independent
of the choice of gluing of $S$. Therefore we have proved that there
is a well-defined map $\bd\co\CW \rightarrow \R^{2c}$
 from a normal class in $\CW$ to a boundary homology class in the
first homology of the boundary tori. As there are assumed to be
$c$ boundary tori, the image is in $\R^{2c}$ as claimed.

It is possible to have a situation where this boundary homology class is
zero but the surface $S$ is spun normal, for example if there are two
essential parallel curves $C, C'$  in the boundary of a core of $S$ and
the surface is spinning around a boundary torus in opposite directions
 from $C, C'$. It is easy to see that this cannot happen for an
embedded surface $S$, but in the general case of singular surfaces, such
behaviour cannot be ruled out. However we can cut off the two infinite
half-open annuli of triangular disks in $S$ and replace by a compact
annulus of triangular disks joining up $C, C'$ to form a normal
surface, which is no longer spun. This illustrates the general fact that
the kernel of $\bd$ consists of classes which can be represented by
standard normal surfaces.

To prove this, if there are several curves $C_1, ...,C_m$ in the
boundary of a core of $S$ and the sum of the homology classes of these
curves is zero, then we can fill in the curves by a cycle consisting of
triangles in the boundary tori of $M$. This may be a singular choice,
since the curves might intersect. The cancelling of the signs corresponds
to the normality of the resulting possibly singular closed surface.

In the non-orientable case, we proceed in a similar manner to the
previous section. Let $\ti M$ be the orientable double covering
with covering transformation $\s$ and lifted triangulation $\ti
\G$. Any spun normal surface lifts to a $\s$--invariant surface or
a pair of such surfaces in $\ti M$. Now the boundary map $\ti \bd$
takes such a surface or pair of surfaces to a $\s$--invariant
homology class in the boundary tori of $\ti M$. Hence we can
project this to a corresponding class in the boundary tori and
Klein bottles of $M$, and hence define the boundary map in $M$.
For a pair of tori interchanged by $\s$ this is easy, since we
have a pair of homology classes in the tori which are the two
lifts of a single class in a boundary torus of $M$. For a single
$\s$--invariant torus covering a Klein bottle, a $\s$--invariant
homology class or pair of classes interchanged by $\s$, both
project to a homology class in the Klein bottle. Notice that the
boundary of a spun normal surface is always a 2--sided curve,
since it is the end of a half open annulus with infinitely many
triangles projecting to a torus or Klein bottle. So we can again
use $\R$ coefficients for convenience, without losing any
information.

The argument in the introduction of this section shows that this map
is onto and so using $\Z$ coefficients, the corresponding map has an
image of finite index in $\Z^{2c}$.

An elegant alternative approach was suggested by the referee. When
considering all possible gluings of a surface from the normal disk
types, it is necessary to observe that the resulting boundary of
the surface is independent of the choice of gluing. However, one
can consider the sum of the oriented boundary arcs of each
quadrilateral as forming a 1--chain. Then the sum of these
quadrilateral boundary arcs is the boundary of the sum of the
quadrilaterals, viewed as a 2--chain, since a pair of oppositely
oriented boundary arcs, coming from two quadrilaterals which may
be glued along these arcs or may not, depending on the way the
normal surface is formed from the disk types, will cancel out in
this 1--chain. Hence the sum, which is a 1--cycle because of the
matching equations and gluing rule, is independent of the choice
of how the quadrilaterals are glued up.
\end{proof}

In the final section we examine
$\bd$ with $\Z$ coefficients more closely, for the simplest
interesting orientable and
non-orientable examples.

\section{The figure--8 knot complement and the Gieseking manifold}
\label{sec5}

In the previous section, we defined the boundary map $\bd\co \CW
\rightarrow \R^{2c}$ and showed that it has rank $2c$. Here the
figure--8 knot complement is presented as an example for which the
boundary map is not onto, using $\Z$ coefficients. Also we
describe the Gieseking manifold as a non-orientable example.

\begin{figure}[ht!]
\cl{\relabelbox\small
\epsfxsize=4.5in
\epsfbox{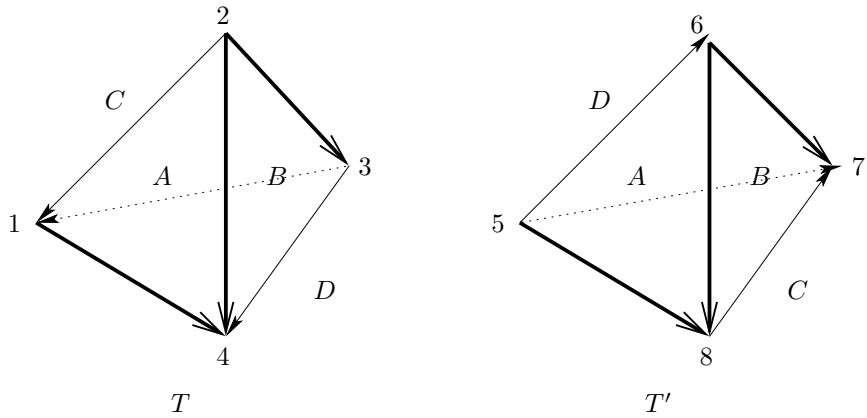}
\relabel{A}{$A$}
\relabel{B}{$B$}
\relabel{C}{$C$}
\relabel{D}{$D$}
\relabel{A'}{$A$}
\relabel{B'}{$B$}
\relabel{C'}{$C$}
\relabel{D'}{$D$}
\relabel{T}{$T$}
\relabel{T'}{$T'$}
\relabel{1}{1}
\relabel{2}{2}
\relabel{3}{3}
\relabel{4}{4}
\relabel{5}{5}
\relabel{6}{6}
\relabel{7}{7}
\relabel{8}{8}
\endrelabelbox}
\caption{An ideal triangulation of the figure--8 knot complement}
\label{fig6}
\end{figure}

\begin{figure}[ht!]
\cl{\relabelbox\small
\epsfxsize=5in
\epsfbox{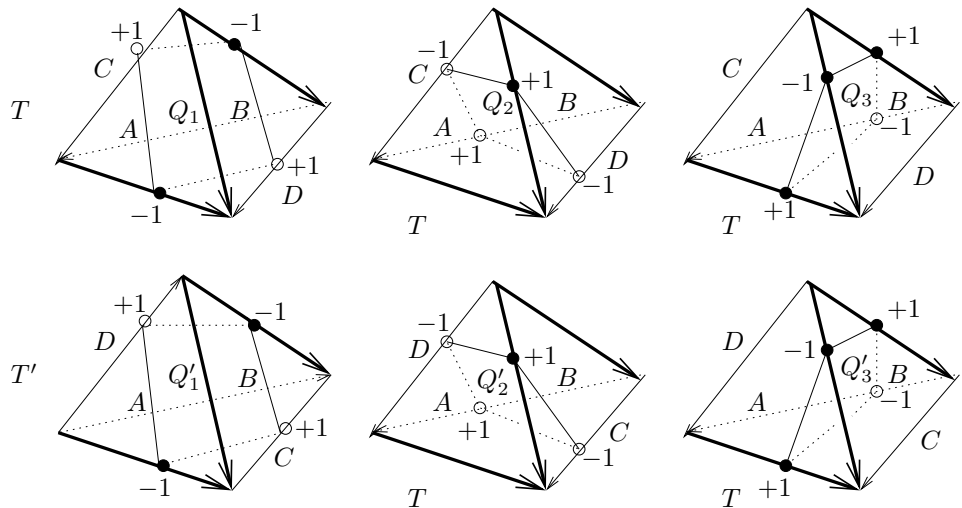}
\adjustrelabel < 0pt,+2pt> {A1}{$A$}
\adjustrelabel <-5pt, 0pt> {B1}{$B$}
\relabel{C1}{$C$}
\adjustrelabel <-4pt,-4pt> {D1}{$D$}
\adjustrelabel <-5pt,+5pt> {-1a}{$-1$}
\adjustrelabel <-2pt,-2pt> {-1b}{$-1$}
\adjustrelabel < 0pt,+2pt> {+1a}{$+1$}
\adjustrelabel <-2pt,-2pt> {+1b}{$+1$}
\adjustrelabel < 0pt,+2pt> {A2}{$A$}
\adjustrelabel <-5pt,+5pt> {B2}{$B$}
\adjustrelabel <-10pt, 0pt> {C2}{$C$}
\relabel{D2}{$D$}
\adjustrelabel <-2pt,-2pt> {-1c}{$-1$}
\adjustrelabel <-2pt,+2pt> {-1d}{$-1$}
\adjustrelabel <-2pt,+2pt> {+1c}{$+1$}
\adjustrelabel <-2pt,+2pt> {+1d}{$+1$}
\adjustrelabel < 0pt,+3pt> {A3}{$A$}
\adjustrelabel < 0pt,+3pt> {B3}{$B$}
\adjustrelabel < 0pt,-2pt> {C3}{$C$}
\adjustrelabel <-3pt,+3pt> {D3}{$D$}
\adjustrelabel <-3pt,+1pt> {-1e}{$-1$}
\adjustrelabel <-6pt,-2pt> {-1f}{$-1$}
\adjustrelabel <-1pt,-1pt> {+1e}{$+1$}
\adjustrelabel <-3pt,-1pt> {+1f}{$+1$}
\adjustrelabel < 0pt,+3pt> {A4}{$A$}
\relabel{B4}{$B$}
\relabel{C4}{$C$}
\adjustrelabel <-3pt,+3pt> {D4}{$D$}
\relabel{-1g}{$-1$}
\relabel{-1h}{$-1$}
\relabel{+1g}{$+1$}
\relabel{+1h}{$+1$}
\adjustrelabel < 0pt,+2pt> {A5}{$A$}
\adjustrelabel < 0pt,+2pt> {B5}{$B$}
\adjustrelabel <-2pt,+3pt> {C5}{$C$}
\relabel{D5}{$D$}
\adjustrelabel <-3pt,+1pt> {-1i}{$-1$}
\adjustrelabel <-5pt, 0pt> {-1j}{$-1$}
\adjustrelabel < 0pt,-2pt> {+1i}{$+1$}
\adjustrelabel < 0pt,-2pt> {+1j}{$+1$}
\adjustrelabel < 0pt,+2pt> {A6}{$A$}
\relabel{B6}{$B$}
\adjustrelabel <-3pt,+3pt> {C6}{$C$}
\relabel{D6}{$D$}
\relabel{-1k}{$-1$}
\relabel{-1l}{$-1$}
\relabel{+1k}{$+1$}
\adjustrelabel <-2pt,-2pt> {+1l}{$+1$}
\relabel{Q1}{$Q_1$}
\relabel{Q2}{$Q_2$}
\adjustrelabel <-2pt, 0pt> {Q3}{$Q_3$}
\relabel{Q'1}{$Q'_1$}
\adjustrelabel <-2pt,-2pt> {Q'2}{$Q'_2$}
\adjustrelabel <-2pt, 0pt> {Q'3}{$Q'_3$}
\relabel{T1}{$T$}
\relabel{T2}{$T$}
\relabel{T3}{$T$}
\relabel{T4}{$T$}
\relabel{T5}{$T$}
\relabel{T'1}{$T'$}
\endrelabelbox}
\caption{All quadrilateral disk types}
\label{fig7}
\end{figure}

Figure \ref{fig6} shows an ideal triangulation of the figure--8 knot
complement described by Thurston \cite{th}. It has two tetrahedra,
two edges and one torus cusp. Following the argument in the
previous section, the vector space $\CW$ of spun and ordinary
normal surfaces in the figure--8 knot complement has dimension 5.
Define the boundary map $\bd\co \CW \rightarrow \R^{2}$ described in
section \ref{sec4}. We compute the total homology classes of some normal
surfaces which imply that the boundary map is not onto. For each
tetrahedron in Figure \ref{fig6}, we have three types of quadrilaterals
denoted by $Q_1$, $Q_2$ and $Q_3$ in $T$ and $Q'_{1}$, $Q'_{2}$
and $Q'_{3}$ in $T'$ (see Figure \ref{fig7}). Using the right hand rule, we
give signs of corners of the quadrilaterals as in Figure \ref{fig7}.

\noindent Let $x_i$ and $y_i$ denote the number of quadrilaterals
of $Q_i$ and $Q'_{i}$ types respectively, for $i=1, 2, 3$. The
following is the system of $Q$--matching equations of the figure--8
knot complement \cite{ka1}, for the two edges of the
triangulation. As in Section \ref{sec3}, the sum of these two equations is
zero, ie one of the two equations is redundant.

\[ \left\{ \begin{array}{l}
               -2x_{1}+x_{2}+x_{3}-2y_{1}+y_{2}+y_{3}=0 \\
               2x_{1}-x_{2}-x_{3}+2y_{1}-y_{2}-y_{3}=0.
             \end{array}
     \right. \]
\noindent There are twenty (possible) fundamental solutions of the
system. Among those, the following six are all the different
solutions up to symmetry.
\[ \begin{array}{l}
       s_{1}=(1,0,0,0,0,2),  \\
       s_{2}=(1,2,0,0,0,0), \\
       s_{3}=(1,1,1,0,0,0), \\
       s_{4}=(1,1,0,0,1,0), \\
       s_{5}=(1,1,0,0,0,1), \\
       s_{6}=(1,0,0,0,1,1). \\
     \end{array}
\]
It is known that the figure--8 knot complement has the symmetry
group isomorphic to the dihedral group
$$D_{4} = \langle f, h | f^{4} = h^{2} = 1, h \circ f \circ h^{-1} =
  f^{-1}\rangle.$$
\noindent So all the other fundamental solutions can be obtained
by isometries in $D_4$. We will show that all the normal surfaces
corresponding to the above six solutions have boundary with an
even number at the first component of the total homology class.
This will imply that the normal surfaces corresponding to all the
other fundamental solutions have also such boundary curves. Hence
the boundary map has image with only even numbers at the first
component, so is not onto. (Our convention will be that the
meridian curve will have homology class $(1,0)$ and the longitude
will correspond to $(0,1)$ for suitable orientations, as described
below).

We label the vertices of the ideal tetrahedra of the figure--8 knot
complement as shown in Figure \ref{fig6}. Then the symmetries $f$ and $h$
can be represented by the following permutations of the eight
vertices;
$$f=(1234)(5876)\qua \text{and}\qua h=(15)(26)(37)(48)$$

\begin{figure}[ht!]
\cl{\relabelbox\small
\epsfxsize=4.5in
\epsfbox{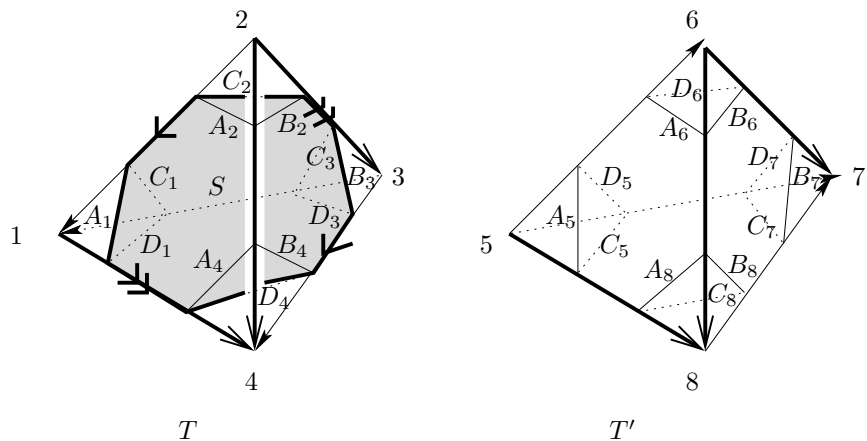}
\relabel{1}{1}
\relabel{2}{2}
\relabel{3}{3}
\relabel{4}{4}
\relabel{5}{5}
\relabel{6}{6}
\relabel{7}{7}
\relabel{8}{8}
\relabel{T}{$T$}
\relabel{T'}{$T'$}
\adjustrelabel <-3pt, 0pt> {A1}{$A_1$}
\adjustrelabel < 0pt,-3pt> {A2}{$A_2$}
\adjustrelabel <-2pt,+2pt> {A3}{$A_4$}
\adjustrelabel < 0pt,-2pt> {B2}{$B_2$}
\relabel{B3}{$B_3$}
\relabel{B1}{$B_4$}
\relabel{C3}{$C_1$}
\adjustrelabel <-2pt,+3pt> {C1}{$C_2$}
\relabel{C2}{$C_3$}
\relabel{D3}{$D_1$}
\adjustrelabel <-3pt, 0pt> {D2}{$D_3$}
\adjustrelabel < 0pt,-1pt> {D1}{$D_4$}
\relabel{S}{$S$}
\adjustrelabel <-2pt, 0pt> {A4}{$A_5$}
\relabel{A5}{$A_6$}
\adjustrelabel <-2pt,+2pt> {A6}{$A_8$}
\relabel{B4}{$B_6$}
\adjustrelabel <-2pt, 0pt> {B5}{$B_7$}
\relabel{B6}{$B_8$}
\relabel{C6}{$C_5$}
\relabel{C4}{$C_7$}
\relabel{C5}{$C_8$}
\relabel{D4}{$D_5$}
\adjustrelabel <-3pt, 0pt> {D6}{$D_6$}
\relabel{D5}{$D_7$}
\endrelabelbox}
\caption{All arc types}
\label{fig8}
\end{figure}

\noindent Table \ref{tbl9} shows the transformations of arc
types in Figure \ref{fig8} by the symmetries $f$ and $h$. The letter label
of each arc type is given by the label of the face containing the
arc type with subscript label being the vertex cut off by the arc
type.

\begin{table}[ht!]
\centering
\begin{tabular}{c|c|c|c|c|c|c|c|c|c|c|c|c}
    & $A_1$ & $A_2$ & $A_4$ & $B_2$ & $B_3$ & $B_4$
    & $C_1$ & $C_2$ & $C_3$ & $D_1$ & $D_3$ & $D_4$ \\ \hline
$f$ & $C_2$ & $C_3$ & $C_1$ & $D_3$ & $D_4$ & $D_1$
    & $B_2$ & $B_3$ & $B_4$ & $A_2$ & $A_4$ & $A_1$ \\ \hline
$h$ & $A_5$ & $A_6$ & $A_8$ & $B_6$ & $B_7$ & $B_8$
    & $D_5$ & $D_6$ & $D_7$ & $C_5$ & $C_7$ & $C_8$
\end{tabular}
\vskip 5mm
\begin{tabular}{c|c|c|c|c|c|c|c|c|c|c|c|c}
    & $A_5$ & $A_6$ & $A_8$ & $B_6$ & $B_7$ & $B_8$
    & $C_5$ & $C_7$ & $C_8$ & $D_5$ & $D_6$ & $D_7$ \\ \hline
$f$ & $C_8$ & $C_5$ & $C_7$ & $D_5$ & $D_6$ & $D_7$
    & $B_8$ & $B_6$ & $B_7$ & $A_8$ & $A_5$ & $A_6$ \\ \hline
$h$ & $A_1$ & $A_2$ & $A_4$ & $B_2$ & $B_3$ & $B_4$
    & $D_1$ & $D_3$ & $D_4$ & $C_1$ & $C_2$ & $C_3$
\end{tabular}

\caption{Transformations of arc types by $f$ and $h$}
\label{tbl9}
\end{table}

\begin{figure}[ht!] 
\centering
\epsfxsize=5in
\epsfbox{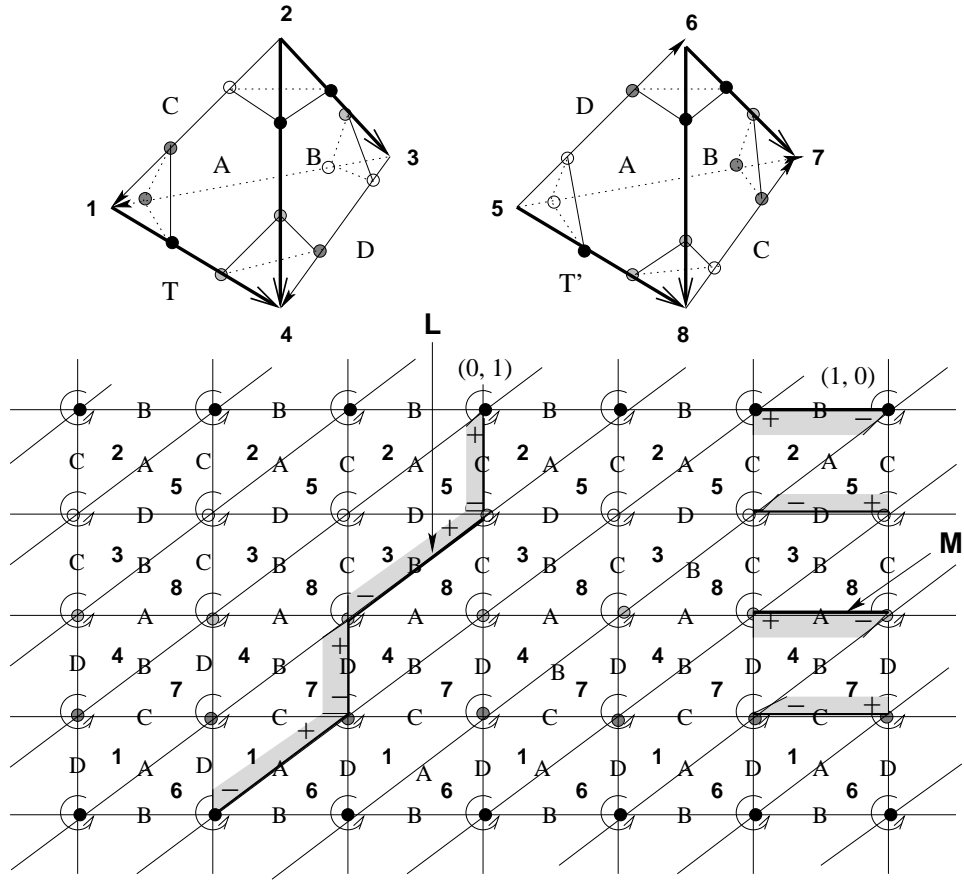}
\caption{Boundary torus}
\label{fig10}
\end{figure}

Let $S$ be the spanning once-punctured torus for the knot which has
the boundary curve $A_1C_2B_3D_4$ as shown in Figure \ref{fig10}. We
let, for convenience, $L=A_1C_5B_3D_7$ be the curve going
(homotopically) around the longitude of the torus cusp once so that it
has the homology class $(0,1)$ (see Figure \ref{fig10}). Notice that
the orientation of the boundary arcs of each quadrilateral is given in
section 4 so that the orientation of the boundary curves of a
quadrilateral surface $S_Q$ is also fixed automatically. From the
table, we can see that $L$ is transformed to the curve
$C_{2}B_{8}D_{4}A_{6}$, which also goes (homotopically) around the
longitude once, by the symmetry $f$. By the symmetry $h$, $L$ is
transformed to $A_{5}D_{1}B_{7}C_{3}$, which traverses the longitude
(homotopically) in the opposite direction to $L$. $M=A_{4}$ is a curve
going (homotopically) around the meridian with homology class $(1,0)$
and it is transformed to $C_1$ by $f$, and to $A_8$ by $h$. Hence the
homology class of the meridian is preserved by both symmetries, but we
need to determine whether its orientation is preserved or reversed.
To do this, consider the curve $D_1D_4C_3C_1$ which follows
(homotopically) the longitude once in the negative orientation and the
meridian 2 times in the negative orientation (see Figure \ref{fig10}),
This curve is transformed to $A_2A_1B_4B_3$ by $f$ and to
$C_5C_8D_7D_6$, by $h$. In the first case, we see that the image curve
winds (homotopically) around the longitude once in the negative
orientation and the meridian twice in the positive orientation,
whereas in the second case the curve is taken to itself
(homotopically) with opposite orientation. So we conclude that $f$
preserves the homology class of $L$ but reverses that of $M$, whereas
$h$ reverses both classes. Hence both $C_1$ and $A_8$ traverse the
meridian (homotopically) in the opposite direction to $A_4$ (see
Figure \ref{fig10}) and have the homology class $(-1,0)$. Also since
$C_1$ is transformed to $B_2$ by $f$, and to $D_5$ by $h$, the
homology classes of $B_2$ and $D_5$ are $(1,0)$. It follows
automatically that $B_6$ and $D_3$ have the homology class
$(-1,0)$. All the other non-trivial curves are linear combinations of
the meridian and longitude.

Now we will compute the boundary curves of the spun and ordinary
normal surfaces corresponding to the solutions $s_{1}, s_{2},
s_{3}, s_{4}, s_{5}$ and $s_{6}$. Since the total homology class
of the boundary curves of any of the surfaces corresponding to a
solution of the matching equations is independent of the choice of
gluing, we choose our surfaces for convenience. Figure \ref{fig10} is
useful as a reference, when computing the boundary class. Notice
that each of these solutions has Euler characteristic $-1$, since
a quadrilateral has `negative curvature' $-2\pi/3$, coming from
its four vertices all lying on edges of degree $6$. Hence we
expect each of the surfaces to be either a once-punctured Klein
bottle, a twice-punctured projective plane or a thrice-punctured
sphere, if all the boundary curves are essential. We can cap off
any inessential boundary curves by disks in the torus at the cusp.

\begin{figure}[ht!]
\centering
\epsfxsize=4.5in
\epsfbox{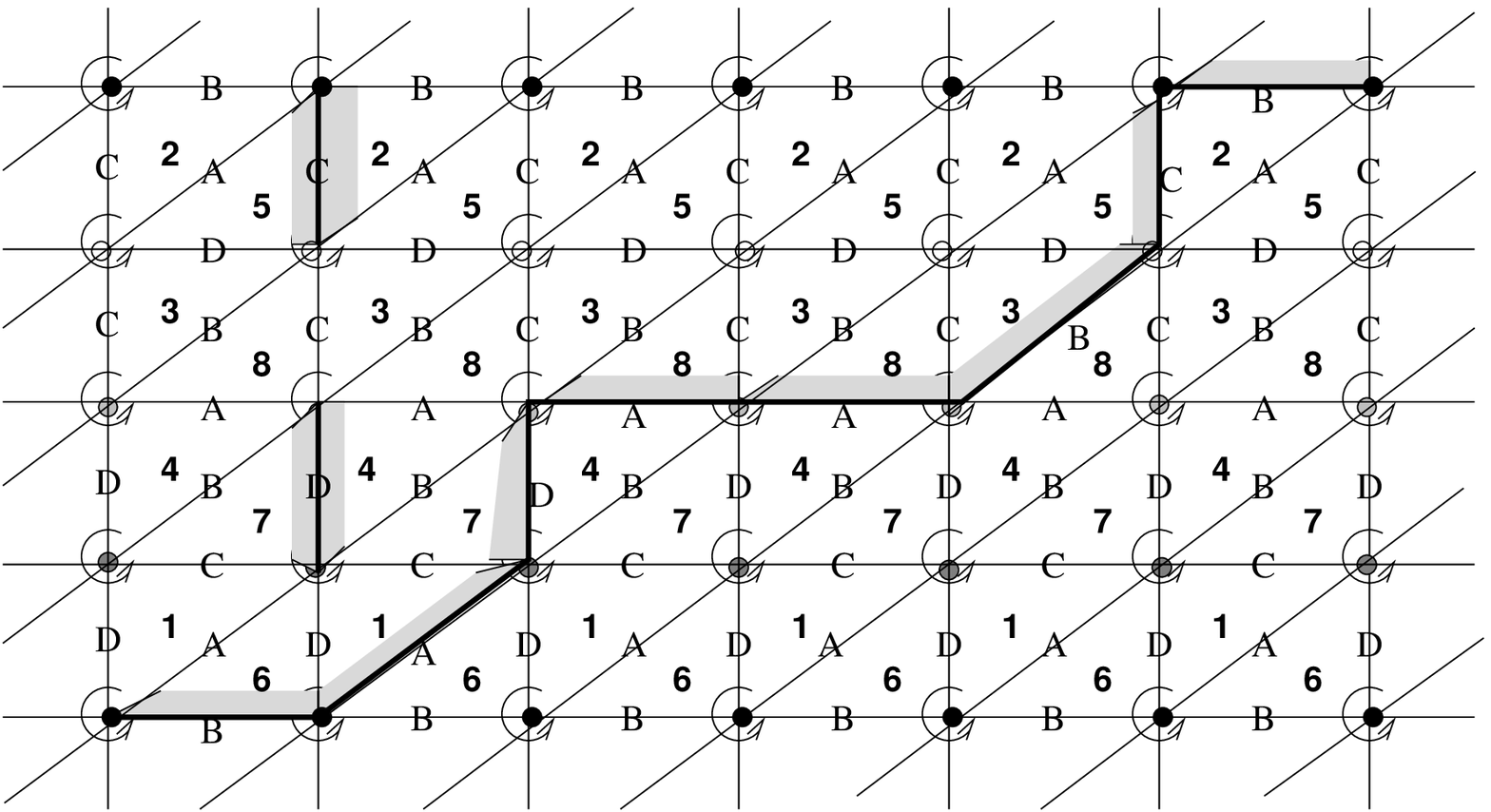}
\caption{Boundary curves of $F_1$}
\label{fig14}
\end{figure}

Let $F_1$ be a normal surface obtained from the solution $s_1$.
The edges of the quadrilaterals of $F_1$ are $A_1$, $C_2$, $B_3$,
$D_4$, $2B_6$, $2A_8$, $2C_5$ and $2D_7$. We may assume that the
edges $C_2$ and $C_5$, and $D_4$ and $D_7$ are glued together and
cancelled out in the boundary class. Hence the boundary of $F_1$
is the curve $A_8D_7A_1B_6B_6C_5B_3A_8$ which has
the total homology class $(-4,1)$, ie winds once around the
longitude and minus four times around the meridian (see Figure \ref{fig14}).
Notice that as is well-known, at this slope, there is an embedded
once-punctured Klein bottle, which is just the surface $F_1$.

\begin{figure}[ht!]
\centering
\epsfxsize=4.5in
\epsfbox{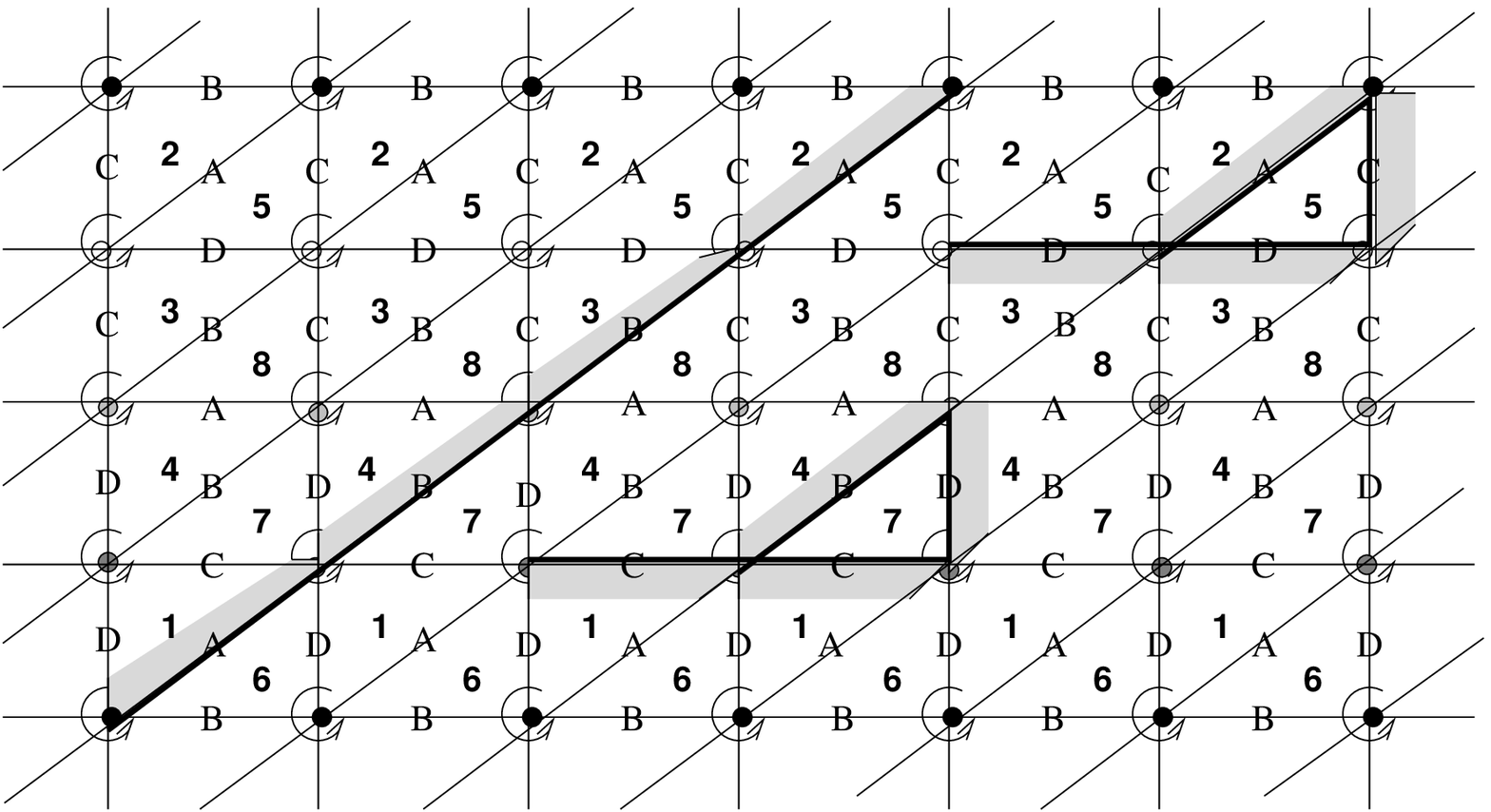}
\caption{Boundary curves of $F_2$}
\label{fig15}
\end{figure}

Now for the solution
$s_2 = (1,2,0,0,0,0)$, a normal surface $F_2$ corresponding to
$s_2$ has the following three boundary components; $A_2B_3B_4A_1$,
$C_1C_1D_4B_4$, and $D_{3}D_{3}C_{2}A_{2}$ which have the homology
classes $(-2,1)$, $(-1,0)$ and $(-1,0)$, respectively (see Figure
\ref{fig15}). Hence the total homology class of the boundary curves of
$F_2$ is $(-4,1)$. $F_2$ must be an immersed thrice-punctured
sphere.

\begin{figure}[ht!] 
\centering
\epsfxsize=4.5in
\epsfbox{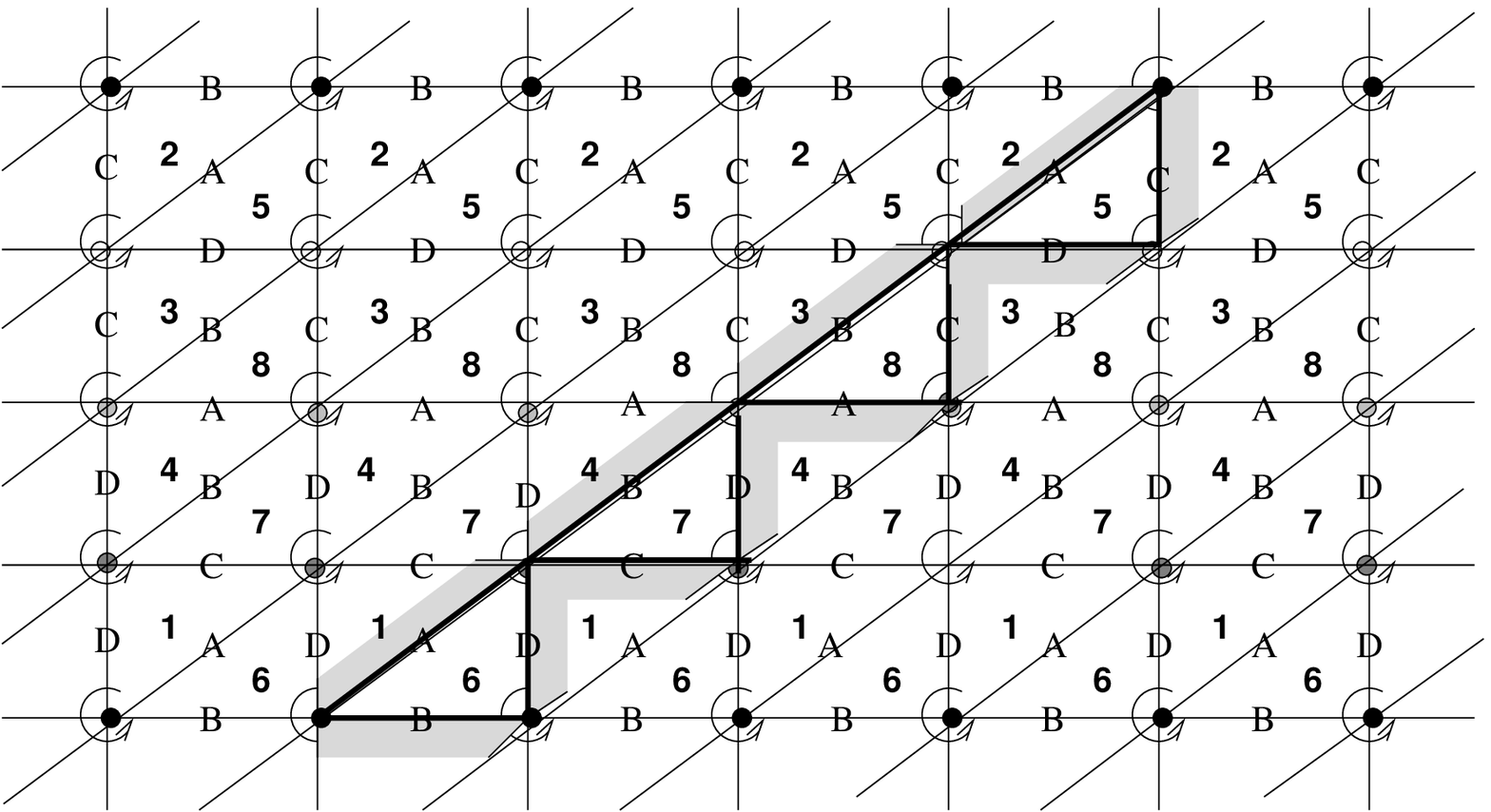}
\caption{Boundary curves of $F_3$}
\label{fig16}
\end{figure}

Let $F_3$ be a normal surface corresponding to the solution $s_3$.  We
can choose the gluing so that $F_3$ has the boundary curves
$A_1D_1B_2$, $B_4D_4C_1$, $B_3C_3A_4$ and $A_2C_2D_3$, all of which
have the homology class (0,0) (see Figure \ref{fig16}).  After capping
off the trivial curves, this is a tetrahedral surface of the canonical
basis defined in the section \ref{sec2}. Notice that the surface has branch
points. Actually there is no gluing for an immersed normal surface
corresponding to the solution $s_3$.

\begin{figure}[ht!]
\centering
\epsfxsize=4.5in
\epsfbox{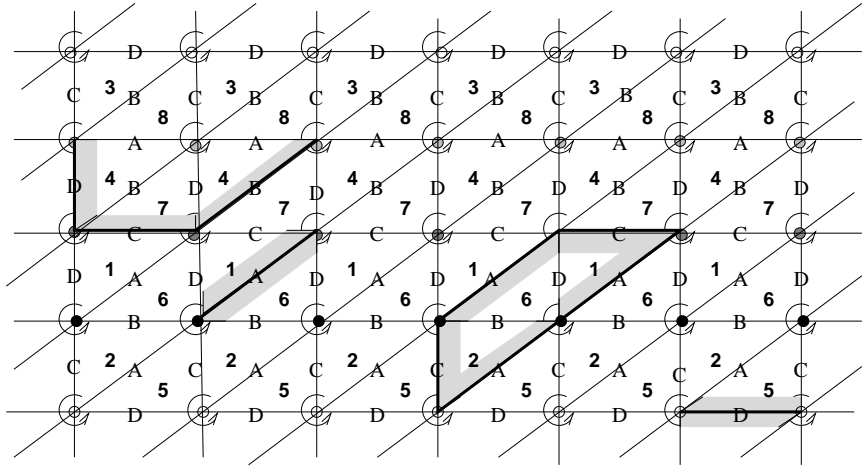}
\caption{Boundary curves of $F_4$}
\label{fig17}
\end{figure}

Let $F_4$ be a normal surface corresponding to the solution $s_4 =
(1,1,0,0,1,0)$. We choose $F_4$ to have the boundary curves
$B_4C_1D_4$ and $C_2A_6C_7A_1A_2$ which have the homology classes
$(0,0)$ and $(2,0)$, respectively (see Figure \ref{fig17}).
Hence the total homology class is $(2,0)$. As the case of $s_3$,
there is no corresponding immersed normal surface to the solution $s_4$.

\begin{figure}[ht!]
\centering
\epsfxsize=4.5in
\epsfbox{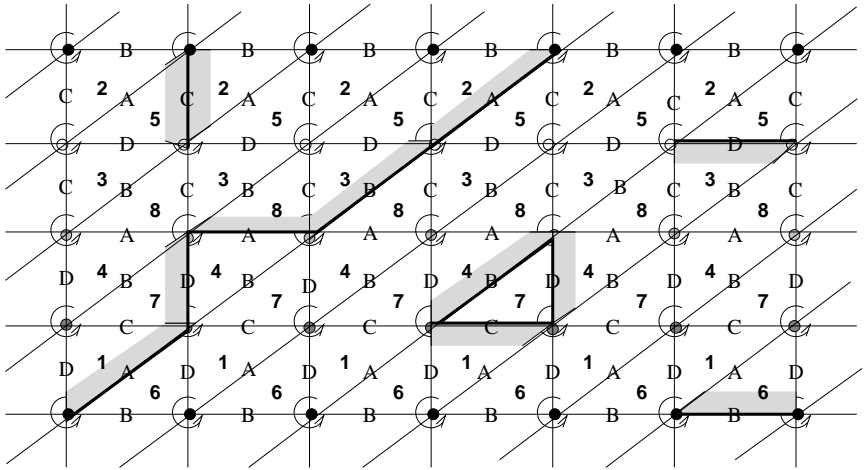}
\caption{Boundary curves of $F_5$}
\label{fig18}
\end{figure}

For the solution $s_5 = (1,1,0,0,0,1)$, a corresponding normal
surface $F_5$ has the boundary curves $A_{2}B_{3}A_{8}D_{7}A_{1}$,
$B_{4}D_{4}C_{1}$, $B_6$ and $D_3$. The homology classes of the
curves are $(-2,1)$, $(0,0)$, $(-1,0)$ and $(-1,0)$ and the total
homology class is $(-4,1)$ (see Figure \ref{fig18}). After capping off the
trivial curve, this is an immersed thrice-punctured sphere.

\begin{figure}[ht!]
\centering
\epsfxsize=4.5in
\epsfbox{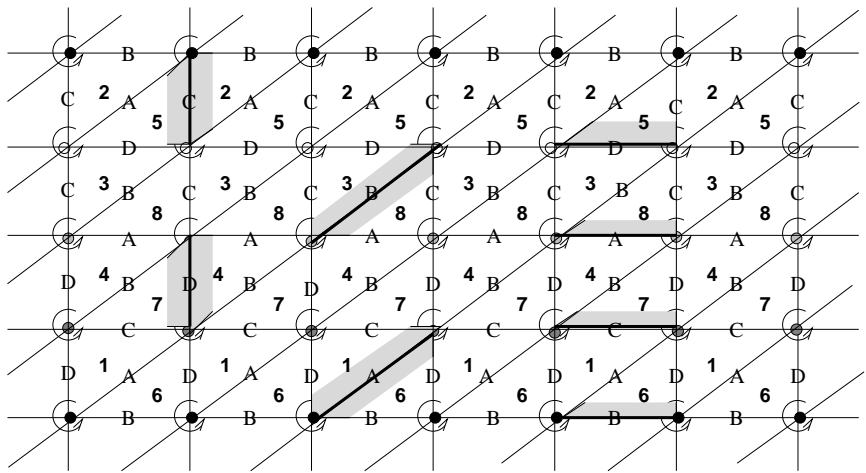}
\caption{Boundary curves of $F_6$}
\label{fig19}
\end{figure}

Finally for the solution $s_6 = (1,0,0,0,1,1)$, we can choose a
normal surface $F_6$ which has boundary curves $A_{8}$, $B_{6}$,
$C_{7}$ and $D_{5}$. Then the homology classes of the curves are
$(-1,0)$, $(-1,0)$, $(1,0)$ and $(1,0)$, respectively and the
total homology class of the boundary curves of $F_6$ is $(0,0)$
(see Figure \ref{fig19}). In this case, the surface has a branch point.
Notice that there is no gluing for an immersed normal surface
corresponding to the solution $s_6$. 

Therefore all six solutions correspond to normal surfaces which have
even numbers at the first component of the total homology class of the
boundary curves.  Since the normal surfaces corresponding to the other
fundamental solutions are obtained from the normal surfaces
corresponding to these six solutions by symmetries, and the symmetries
transform the boundary curves as shown previously, the total homology
classes of the boundary curves of any spun normal surfaces in the
figure--8 knot complement have only even numbers at the first
component. Hence the boundary map $\bd\co \CW \rightarrow \R^{2}$ has
image the subgroup of index $2$ consisting of all slopes of the form
$(2m,n)$.

\begin{figure}[ht!]
\cl{\relabelbox\small
\epsfxsize=4.2in
\epsfbox{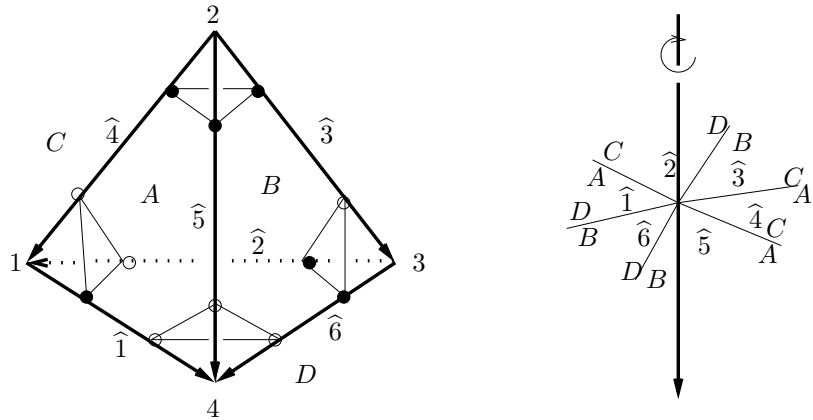}
\relabel{1a}{1}
\relabel{2a}{2}
\relabel{3a}{3}
\relabel{4a}{4}
\relabel{1b}{$\widehat{1}$}
\relabel{2b}{$\widehat{2}$}
\relabel{3b}{$\widehat{3}$}
\relabel{4b}{$\widehat{4}$}
\adjustrelabel <-3pt, 0pt> {5b}{$\widehat{5}$}
\relabel{6b}{$\widehat{6}$}
\relabel{1c}{$\widehat{1}$}
\relabel{2c}{$\widehat{2}$}
\relabel{3c}{$\widehat{3}$}
\relabel{4c}{$\widehat{4}$}
\relabel{5c}{$\widehat{5}$}
\relabel{6c}{$\widehat{6}$}
\relabel{A1}{$A$}
\relabel{A2}{$A$}
\relabel{A3}{$A$}
\relabel{A4}{$A$}
\relabel{B1}{$B$}
\relabel{B2}{$B$}
\relabel{B3}{$B$}
\relabel{B4}{$B$}
\relabel{C1}{$C$}
\relabel{C2}{$C$}
\relabel{C3}{$C$}
\relabel{C4}{$C$}
\relabel{D1}{$D$}
\relabel{D2}{$D$}
\relabel{D3}{$D$}
\relabel{D4}{$D$}
\endrelabelbox}
\caption{Gieseking manifold}
\label{fig11}
\end{figure}

Finally we study the Gieseking manifold, which has an ideal
triangulation with one tetrahedron, one edge and a Klein bottle
cusp. The Gieseking manifold is double-covered by the figure--8
knot complement and has covering transformation by the symmetry $f
\circ h = (18)(25)(36)(47)$. Figure \ref{fig11} shows the Gieseking
manifold and a fixed orientation going around the unique edge.

\begin{figure}[ht!]
\cl{\relabelbox\small
\epsfxsize=5.2in
\epsfbox{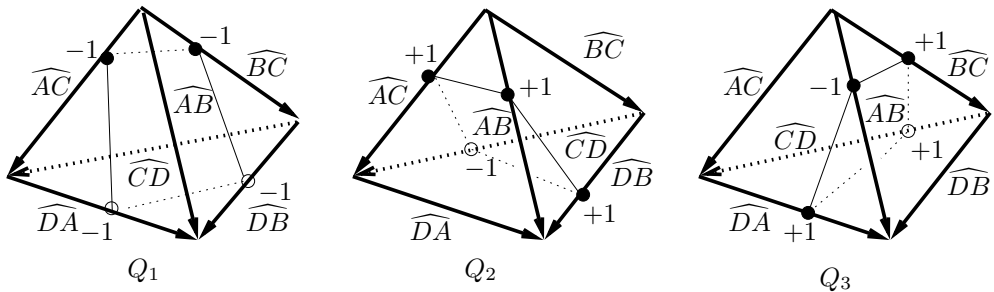}
\adjustrelabel <-4pt,+2pt> {-1a}{$-1$}
\adjustrelabel <-4pt, 0pt> {-1b}{$-1$}
\adjustrelabel <-4pt, 0pt> {-1c}{$-1$}
\adjustrelabel <-3pt, 0pt> {-1d}{$-1$}
\adjustrelabel <-2pt,-2pt> {-1e}{$-1$}
\relabel{-1f}{$-1$}
\relabel{+1a}{$+1$}
\adjustrelabel <-2pt, 0pt> {+1b}{$+1$}
\relabel{+1c}{$+1$}
\relabel{+1d}{$+1$}
\adjustrelabel <-2pt, 2pt> {+1e}{$+1$}
\relabel{+1f}{$+1$}
\relabel{Q1}{$Q_1$}
\relabel{Q2}{$Q_2$}
\relabel{Q3}{$Q_3$}
\adjustrelabel <-2pt, 0pt> {AC1}{$\widehat{AC}$}
\adjustrelabel <-2pt, 0pt> {AC2}{$\widehat{AC}$}
\adjustrelabel <-2pt, 0pt> {AC3}{$\widehat{AC}$}
\relabel{CD1}{$\widehat{DA}$}
\relabel{CD2}{$\widehat{DA}$}
\relabel{CD3}{$\widehat{CD}$}
\relabel{CD4}{$\widehat{DA}$}
\relabel{CD5}{$\widehat{CD}$}
\relabel{CD6}{$\widehat{CD}$}
\relabel{BC1}{$\widehat{BC}$}
\relabel{BC2}{$\widehat{BC}$}
\relabel{BC3}{$\widehat{BC}$}
\relabel{DB1}{$\widehat{DB}$}
\relabel{DB2}{$\widehat{DB}$}
\relabel{DB3}{$\widehat{DB}$}
\relabel{AB1}{$\widehat{AB}$}
\adjustrelabel <-4pt,+2pt> {AB2}{$\widehat{AB}$}
\relabel{AB3}{$\widehat{AB}$}
\endrelabelbox}
\caption{Signs of corners of quadrilaterals in the Gieseking manifold}
\label{fig12}
\end{figure}

We denote each wedge $\widehat{1}, \widehat{2}, \widehat{3}, \widehat{4},
\widehat{5}$
and $\widehat{6}$ by $\widehat{DA}, \widehat{CD}, \widehat{BC}, \widehat{AC},
\widehat{AB}$ and $\widehat{DB}$ according to the labels of the two
adjacent faces and the fixed orientation going around the edge.
Since the wedges $\widehat{CD}, \widehat{AC}$ and $\widehat{DB}$ have the
opposite orientations to the ones obtained by the right hand rule
applied to each edge in the tetrahedron without identification
(see the case of the figure--8 knot complement), the signs of
corners of quadrilaterals are as shown in Figure \ref{fig12}. Therefore
the $\Q$--matching equation of the Gieseking manifold is
$$\begin{array}{l} -2x_{1}+x_{2}+x_{3}=0, \end{array}$$
where $x_1$, $x_2$ and $x_3$ are the number of quadrilaterals of
type $Q_1$, $Q_2$ and $Q_3$, respectively.
Then all positive solutions are generated by the following solutions
$$\begin{array}{l}
       t_{1}=(1,1,1), \\
       t_{2}=(1,2,0), \\
       t_{3}=(1,0,2). \\
     \end{array}$$
Now we will compute the homology class of the boundary curves of
normal surfaces corresponding to each solution $t_1$, $t_2$ and
$t_3$. We will follow exactly the same procedure as in the case of
the figure--8 knot complement.

\begin{figure}[ht!]
\centering
\epsfxsize=5in
\epsfbox{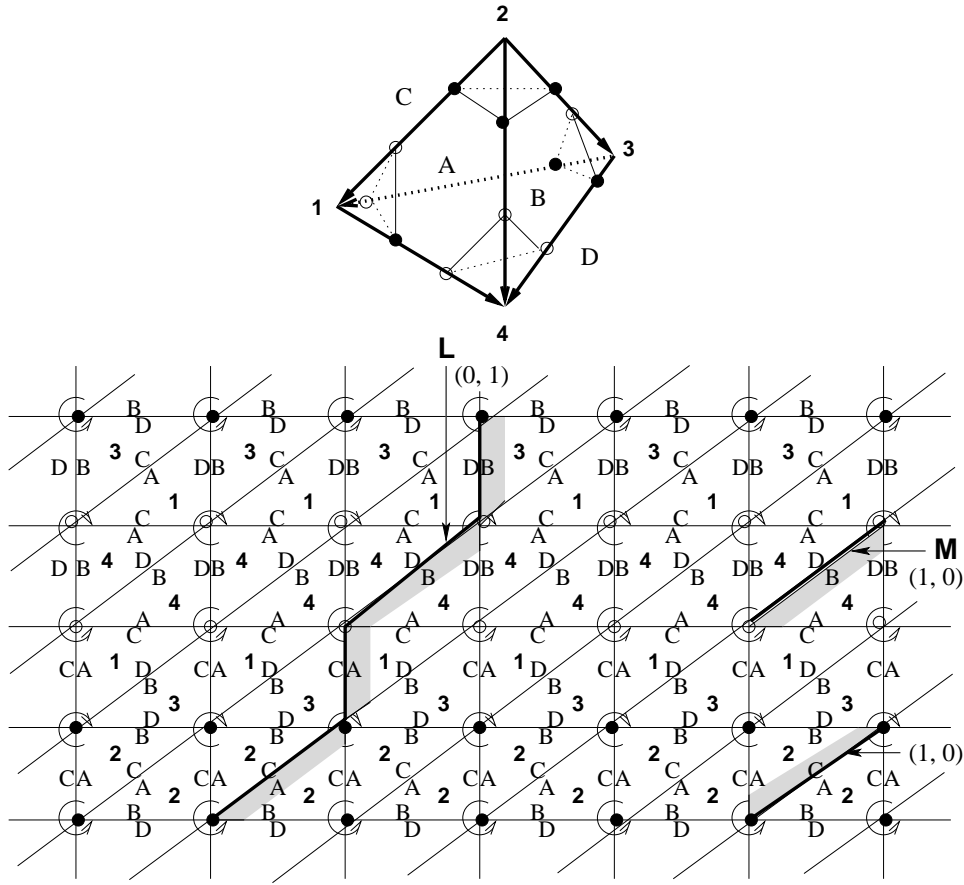}
\caption{Boundary Klein bottle of the Gieseking manifold}
\label{fig13}
\end{figure}

\noindent Figure \ref{fig13} shows the Klein bottle cusp of the Gieseking
manifold. Note that the orientation around each vertex along a
column is shown in opposite directions alternately. This will
assist in determining the signs of the homology classes of curves.
Also we call the longitude the slope $L=A_2A_1B_4B_3$, which is
equivalent to the curve $A_{1}D_{4}B_{3}C_{2}$, and the meridian
the slope $M=B_4$ to conform to the figure--8 knot terminology (see
Figure \ref{fig13}). Notice that $B_4$ lifts to $B_4D_4$, which is
homotopic to the meridian slope $A_4$ of the torus cusp in the
figure--8 knot complement, by the double covering map from the
figure--8 knot complement onto the Gieseking manifold. Let the
homology classes of $L$ and $M$ be $(0,1)$ and $(1,0)$,
respectively. Then we get the homology classes of the following
curves; $A_2\sim \overline{C_2}\sim D_4\sim \overline{B_4}$ with
$(-1,0)$, $D_3\sim A_4 \sim \bar{D_4}B_4$ with $(2,0)$, $A_4\sim
\overline{C_1}\sim \overline{B_2}$ with $(2,0)$, $C_3B_3\sim C_1\sim
B_2$ with $(-2,0)$ and $A_1D_1\sim A_4\sim D_3$ with $(2,0)$.

\begin{figure}[ht!]
\centering
\epsfxsize=4.5in
\epsfbox{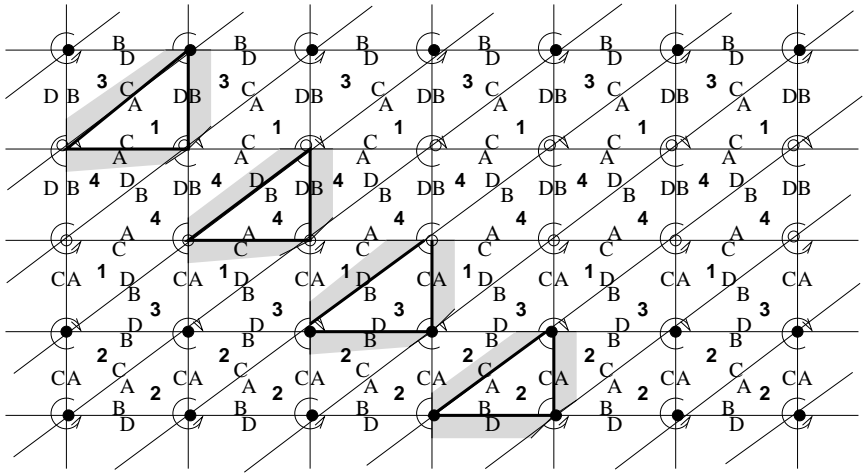}
\caption{Boundary curves of $G_1$}
\label{fig20}
\end{figure}

\begin{figure}[ht!]
\centering
\epsfxsize=4.5in
\epsfbox{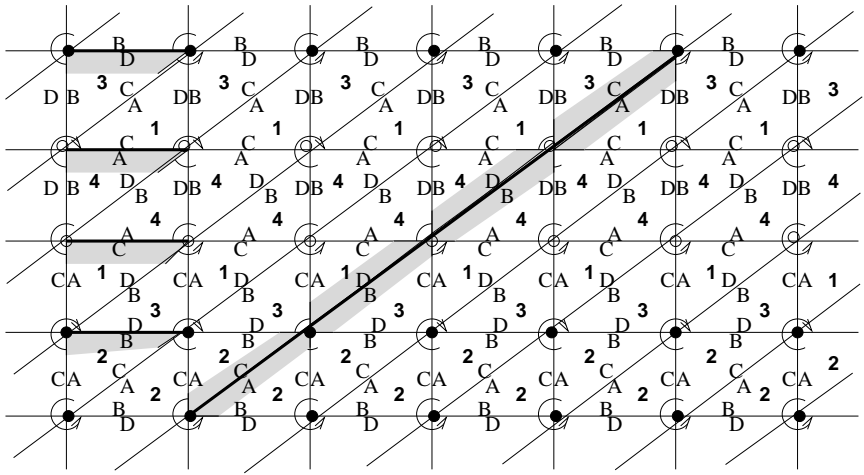}
\caption{Boundary curves of another $G_1$}
\label{fig21}
\end{figure}

Let $G_1$ be a normal surface obtained from the solution $t_1 =
(1,1,1)$. We can choose $G_1$ to have boundary curves
$D_{3}C_{3}B_{3}$, $A_{4}D_{4}B_{4}$, $C_{1}D_{1}A_{1}$ and
$B_{2}C_{2}A_{2}$, all of which have the trivial homology class
(see Figure \ref{fig20}). Hence the total homology class of the boundary
curves of $G_1$ is $(0,0)$. Note that the curves, $D_3$, $A_4$,
$C_1$, $B_2$, $C_{2}D_{1}D_{4}C_{3}$ and $A_{2}A_{1}B_{4}B_{3}$,
are also the boundary curves of a normal surface obtained from the
same solution but by another gluing. If we examine the
orientations of the curves carefully, we will see that the curves
have the homology class $(2,0)$, $(2,0)$, $(-2,0)$, $(-2,0)$,
$(-4,-1)$ and $(4,1)$, respectively (see Figure \ref{fig21}). Hence the
total homology class of these boundary curves is also trivial.

\begin{figure}[ht!]
\centering
\epsfxsize=4.5in
\epsfbox{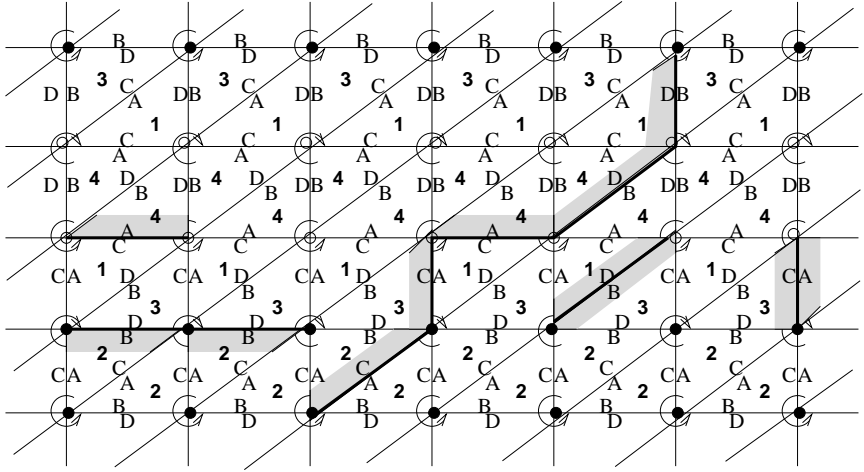}
\caption{Boundary curves of $G_2$}
\label{fig22}
\end{figure}

Next look at the solution $t_2 = (1,0,2)$. We get a normal surface
with boundary curves $C_{2}C_{3}A_{4}D_{4}D_{1}$, $B_{2}B_{2}$ and
$A_{4}$ which have the homology classes $(2,-1)$, $(-4,0)$ and
$(-2,0)$, respectively (see Figure \ref{fig22}). Then the total homology
class is $(0,-1)$. Note that the edge pairs $D_1$ and $B_3$, and $C_3$
and $A_1$ cancel out.

\begin{figure}[ht!]
\centering
\epsfxsize=4.5in
\epsfbox{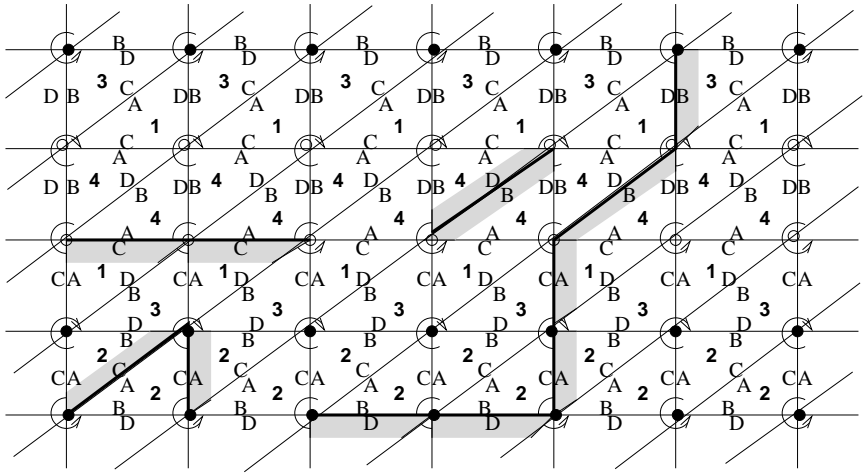}
\caption{Boundary curves of $G_3$}
\label{fig23}
\end{figure}

Finally for the solution $t_3 = (1,2,0)$, we get the following boundary
curves: $D_{3}D_{3}A_{2}A_{1}B_{4}B_{3}$, $C_{2}A_{2}$ and $C_1C_1$
which have the homology class $(6,1)$, $(-2,0)$ and $(-4,0)$ (see
Figure \ref{fig23}). The edge pair $D_4$ and $B_4$ cancels out. Hence the
total homology class is $(0,1)$. Therefore as in the case of the
figure--8 knot complement, the boundary map of the Gieseking
manifold also has only multiples of 4 at the first component and
is not onto. Note that for the first homology of the Klein bottle,
the longitude is an element of order two. Hence we have found the
image of the boundary map is a subgroup of index 4 of the group
$\Z \oplus \Z_2$.

\rk{Acknowledgements}
We would like to thank Stefan Tillmann for helpful
conversations on spun normal surface theory and also an
anonymous referee for a number of very helpful comments, which
have improved the exposition. The second author is supported by a
grant from the Australian Research Council.

\Addresses\recd

\end{document}